\input  amstex
\input amsppt.sty
\magnification1200
\vsize=23.5truecm
\hsize=16.5truecm
\NoBlackBoxes

\def\crp{\overline{\Bbb R}_+}
\def\crm{\overline{\Bbb R}_-}
\def\crpm{\overline{\Bbb R}_\pm}
\def\rn{{\Bbb R}^n}
\def\rnp{{\Bbb R}^n_+}
\def\rnm{\Bbb R^n_-}

\def\crnp{\overline{\Bbb R}^n_+}
\def\crnm{\overline{\Bbb R}^n_-}

\def\comega{\overline\Omega }
\def\ang#1{\langle {#1} \rangle}
\def\rp{ \Bbb R_+}

\def\OP{\operatorname{OP}}
\def\N{\Bbb N}
\def\R{\Bbb R}
\def\C{\Bbb C}

\def\ol{\overline}
\def\E{\Cal E}
\def\F{\Cal F}
\def\simto{\overset\sim\to\rightarrow}
\def\supp{\operatorname{supp}}
\def\chp{\Cal H^+}
\def\chm{\Cal H^-}

\document
\topmatter
\title
Integration by parts for nonsymmetric fractional-order operators on a halfspace 
\endtitle
\author Gerd Grubb \endauthor
\affil
{Department of Mathematical Sciences, Copenhagen University,\\
Universitetsparken 5, DK-2100 Copenhagen, Denmark.\\
E-mail {\tt grubb\@math.ku.dk}}\endaffil
\rightheadtext{integration by parts}
\abstract
For a strongly elliptic pseudodifferential operator $L$ of order $2a$
($0<a<1$) with real kernel, we show an integration-by-parts formula
for solutions of the homogeneous Dirichlet problem, in
the model case where the operator is $x$-independent with homogeneous
symbol, considered on the halfspace $\rnp$. The new aspect compared
to $(-\Delta )^a$ is that $L$ is nonsymmetric,  having
both an even and an odd part. Hence it satisfies a $\mu
$-transmission condition where generally $\mu \ne a$.

We present a complex
method, relying on a factorization in factors holomorphic in $\xi _n$
in the lower or upper complex halfplane, using order-reducing
operators combined with a decomposition
principle originating from Wiener and Hopf. This is in contrast to  a
real, computational method presented very recently by Dipierro,
Ros-Oton, Serra and Valdinoci. Our method allows $\mu $ in a larger
range than they consider.

Another new contribution is the (model) study of ``large'' solutions of
nonhomogeneous Dirichlet problems when $\mu >0$. Here we deduce a ``halfways Green's formula'' for $L$:
$$
\int_{\rnp} Lu\,\bar v\,dx-\int_{\rnp}u\,\overline{  L^*v}\,dx=c\int_{{\Bbb
R}^{n-1}}\gamma _0(u/x_n^{\mu -1 })\,{\gamma _0(\bar v/x_n^{\mu ^*})}\,
dx',
$$ 
when $u$ solves a nonhomogeneous Dirichlet problem for $L$, and $v$
solves a homogeneous Dirichlet problem for $L^*$; $\mu ^*=2a-\mu $.
Finally, we show a full Green's formula, when both $u$ and $v$ solve
nonhomogeneous Dirichlet problems; here both Dirichlet and Neumann
traces of $u$ and $v$ enter, as well as a first-order pseudodifferential operator
over the boundary.
\endabstract



\endtopmatter

\subhead 1. Preface \endsubhead

The fractional Laplacian $(-\Delta )^a$ ($0<a<1$) and its
generalizations have received much attention lately, because
of applications in probability, finance, mathematical physics and
differential geometry. (References to applications are listed e.g.\ in our earlier works \cite{G14-G19}). It is interesting to observe that the operator has been attacked from several
angles: Roughly speaking, 1) there is a probabilistic approach, based on the fact that it
is the infinitesimal generator of a L\'evy process, 2) there is  an approach by
potential theoretic methods searching for similarities with the usual
Laplacian $\Delta $, and 3) there is an approach by Fourier
transformation using that the fractional Laplacian is a pseudodifferential operator. The
first two methods have filled the major part of contributions in recent years; they
deal mainly with {\it real} functions and {\it real} integral
operators (singular integral operators with real kernels). The third
method gave some early results in the 1960s and '70s, but has only been
reintroduced fairly recently for this particular type of operators; it
uses {\it complex} functions and distribution theory, and {\it complex} symbols
and kernels of operators, in an
essential way.
 
In the study of boundary value problems for these operators on open subsets $\Omega $ of $\rn$,
there are often two distinguished steps in the procedure: One step is to
analyze a {\it model case}, where the domain $\Omega $ has the simple
form of a ball or a halfspace, and the operator is
translation-invariant; the other step is to piece the resulting information
together to treat cases where $\Omega $ is more general, or the
operator is $x$-dependent, or both things happen.

The purpose of the present paper is to show how the complex method 
can be used to advantage to  give satisfactory results in the model
case for operators $L$ that generalize the fractional Laplacian by
being nonsymmetric, but still have real homogeneous kernels. (An
example is $L=(-\Delta )^\frac12+\pmb b\cdot \nabla$, $\pmb b\in\rn$;
the fractional Laplacian with critical drift.)

A direct inspiration for this work was a correspondence with the authors
X.\ Ros-Oton, S.\ Dipierro, J.\ Serra and E.\ Valdinoci of the recent
paper \cite{DRSV20}, treating generators $L$ of $\alpha $-stable
processes and
in particular showing an integration-by-parts formula. For the proof of
the formula in the model case of the operator on a halfspace, they use
an entirely real method passing through cumbersome integral operator
formulas and calculations with
special functions, taking up many pages. As we
suggested to them, there is a complex method relying on the knowledge around
pseudodifferential operators, leading to the result in an instructive
way.

The complex method is presented here. The essence of the difference
between the strategy of \cite{DRSV20} and this work
lies in the decomposition of $L$, where \cite{DRSV20} uses a
factorization into two real factors $L=L^\frac12 L^\frac12$, and we decompose $L$ in terms
with symbol  holomorphic in the lower, resp.\ upper complex plane
(by principles originating from  Wiener and Hopf).

An additional gain is that we do not need certain restrictions on the
numerical range of the symbol of $L$ imposed in \cite{DRSV20}, hence
can allow negative transmission numbers $\mu $ or $\mu ^*$, in the
considerations of homogeneous Dirichlet problems.

As a new contribution, we moreover consider ``large'' solutions (with the
singularity $x_n^{\mu -1}$ at the boundary), solving nonhomogeneous
Dirichlet problems. 
Here we show when $\mu ,\mu ^*>0$ how a  ``halfways Green's formula'' can be
deduced from the integration-by-parts formula. Furthermore we work out a full
Green's formula with nontrivial Dirichlet and Neumann traces.

The present note is only concerned with model cases on the halfspace
$\rnp$. For the
transition to cases with more general domains $\Omega $, the
localization method of
\cite{DRSV20} would most likely be a useful tool.

\subsubhead{Outline}\endsubsubhead
In Section 2, the operator $L=\OP(\Cal A(\xi
)+i\Cal B(\xi ))$ is introduced, and a value $\mu $ for which it
satisfies the $\mu $-transmission condition is determined. Section 3
recalls mapping properties and regularity results from \cite{G15} for $L$ on $\rnp$
with homogeneous Dirichlet condition. In Section 4, the
integration-by-parts formula is shown in the model case on $\rnp$ by
an extension of methods from \cite{G16} and \cite{G18}. Section 5
introduces large solutions and the nonhomogeneous Dirichlet problem, and shows how the
result of Section 4 implies a halfways
Green's formula. Section 6 gives the proof of the full Green's formula,
drawing on methods from \cite{G18}. In Section 7, the formulas are extended to
general functions $u,u',v$ in Sobolev-type and H\"older-type function spaces.

\subhead 2. Introduction \endsubhead

In the following we consider pseudodifferential operators $P=\OP(p(\xi ))$ on $\rn$,
where the symbols are independent of $x$:
$$
Pu=\Cal F^{-1}(p(\xi )(\Cal Fu)(\xi )),\text{ where }\Cal F
u= \hat
u(\xi )=
\int_{{\Bbb R}^n}e^{-ix\cdot \xi }u(x)\, dx.\tag2.1
$$
$P$ can also be written as a singular integral operator
$$
(Pu)(x)=PV\int_{\rn}(u(x)-u(x+y))K(y)\, dy, \quad K(y)=\Cal F^{-1}p,\tag2.2
$$
where $PV$ indicates principal value.

In a recent paper \cite{DRSV20}, Dipierro, Ros-Oton, Serra and
Valdinoci have considered such problems without the customary assumption
that $P$ is symmetric; we shall study such cases more closely.

Assume that $p(\xi )$ is $C^\infty $ for $\xi \ne 0$ and homogeneous of degree $m>0$; then $K(y)$ is a
distribution 
homogeneous of degree $-n-m$. \cite{DRSV20} moreover assumes that $K$
is {\bf real} (in order to work with real integral operators), but without
the earlier imposed assumption that $K$ should be {\bf even}, i.e.\
invariant under the mapping $y\mapsto -y$. In the present case, $K$ splits in an
even and an odd part,
$$
K(y)=K_e(y)+K_o(y),\quad K_e(y)=\tfrac12(K(y)+K(-y)),\; K_o(y)=\tfrac12(K(y)-K(-y)).\tag2.3
$$
It is checked from the properties of the Fourier transform that here $p(\xi
)=\Cal A(\xi )+i\Cal B(\xi )$, where $\Cal A(\xi
)=\operatorname{Re}p(\xi )=\Cal F K_e$  and  $\Cal B(\xi
)=\operatorname{Im}p(\xi )=\frac1i\Cal F K_o$ are real, and $\Cal A(\xi )$ is even, $\Cal B(\xi
)$ is odd. The assumptions in
\cite{DRSV20} imply that $p(\xi )$ is strongly elliptic,
i.e., $\Cal A(\xi )>0$ for $\xi \ne 0$. Let us give $p$ the name $\Cal
L(\xi )=\Cal A(\xi )+i\Cal B(\xi )$ and denote the operator
$L=\OP(\Cal L(\xi ))$.

It is moreover assumed in \cite{DSRV20} that $L$ is the infinitesimal
generator of an $\alpha $-stable $n$-dimensional L\'evy process; this
puts certain limitations on the numerical range of $\Cal L$ (cf.\ (2.9)
below) and the
form of the integral operator, which we need not impose in our treatment.
We shall simply assume, in addition to the homogeneity, that  $\Cal L(\xi )$ is  strongly elliptic
and that $K(y)=\Cal F^{-1}\Cal L$ is real.

For a homogeneous symbol $p(\xi )$ of degree $m$, the  $\mu $-transmission
condition with respect to the boundary of $\rnp$ was defined in
\cite{G15, Sect. 3} (with earlier input from H\"ormander
\cite{H66,H85} and Eskin \cite{E81}) as a measurement of how much $p$ deviates from being
symmetric on the normal $(\xi _1,\dots,\xi _{n-1},\xi _n)=(0,1)$ to $\partial\rnp$:
$$
p(0,-1)=e^{i\pi (m-2\mu )}p(0,1)\tag2.4
$$
(note that this only determines $\mu $  mod 1).
When $p$ is {\bf even} in $\xi _n$, (2.4) holds with $\mu =m/2$. The
case $m=2a>0, \mu =a$ has been amply treated in works of Grubb, Ros-Oton
and Serra with coauthors, and many others. In particular, when $p$ is
even in $\xi $, (2.4)
holds with $\mu =m/2$ also for rotations of $\rnp$; hence on a general
domains $\Omega \subset \rn$ it holds regardless of the direction of
the normal to $\partial\Omega $.

What is the index $\mu $ in the case $\Cal L(\xi )=\Cal A(\xi )+i\Cal B(\xi )$?
The order $m$ will be denoted $2a$, $0<a<1$. We have: 
$$
\Cal L(0,1)=\Cal A(0,1)+i\Cal B(0,1),\quad \Cal L(0,-1)=\Cal A(0,-1)+i\Cal
B(0,-1)= \Cal A(0,1)-i\Cal B(0,1),\tag2.5 
$$
in short,
$$
\Cal A(0,1)^{-1}\Cal L(0,\xi _n )=\cases 1+ib\text{ for }\xi _n =1\\
1-ib\text{ for }\xi _n =-1\endcases, \text{ with }b=\frac {\Cal B(0,1)}{\Cal A(0,1)}.\tag2.6
$$
In particular, $\Cal L(0,1)$ and $\Cal L(0,-1)$ have the same length. The angle $\theta $
between the positive real axis and $\Cal L(0,1)$ is
$$
\theta =\operatorname{Arctan} b,\text{ and we set }\pi ^{-1} \theta
=\delta , \text{ whereby }\Cal L(0,1)=e^{i\pi \delta }|\Cal L(0,1)|.\tag2.7
$$
Moreover, 
$$
\Cal L(0,-1)=e^{-i2\theta }\Cal L(0,1)=e^{-i\pi 2\delta }\Cal
L(0,1)=e^{i\pi (2a-2(a+\delta ))}\Cal L(0,1),
$$
so (2.4) holds with 
$$
\mu =a+\delta ;\quad
\delta =\pi ^{-1}\operatorname{Arctan} \frac {\Cal B(0,1)}{\Cal
A(0,1)}\in \,]-\tfrac12,\tfrac12[\, \tag2.8
$$
(all other possible $\mu $ equal this one plus an integer). 
As far as we
understand the analysis in \cite{E81, Ex.\ 6.1}, the value $\mu $ in (2.8) equals the {\it
factorization index}, called $\kappa $ in \cite{E81} and $\mu _0$ in \cite{G15}.

There is the difficulty in comparisons with the Russian works Eskin
\cite{E81} and the subsequent developments by Shargorodsky \cite{S95} and Chkadua and Duduchava
\cite{CD01}, that they use a Fourier transformation having the exponential
factor $e^{+ix\cdot\xi }$ with plus instead of minus. This does not
change the overall objective, but makes it hard to check exact
formulas.

In the notation of \cite{DRSV20}, $m/2$ is called $s$ (or $\alpha $),
and our formula (2.8) for $\mu $
corresponds to their
formula (1.9) for $\gamma $, consistent also with a related formula in
Fernandez-Real and Ros-Oton \cite{FR18, Prop.\ 2.2} concerning the
special case $L=(-\Delta )^\frac12+\pmb b\cdot \nabla$, $\pmb
b\in\rn$, where $a=\frac12$. The extra
restrictions imposed in \cite{DRSV20} require that (with $0<a<1$)
$$
\mu \in \, ]0,2a[\,\cap \,]2a-1,1[\, ,\tag2.9
$$
for their operators. We only assume $a\in \,]0,1[\,$ and $|\delta
|<\frac12$ (whereby $\mu $ can take any value in $\,]-\frac12, \frac32[\,$).

For completely general strongly elliptic symbols $p(\xi )$ (where one drops the
requirement of having a real kernel), $p(0,1)$ and $p(0,-1)$ can
have different lengths, so that $\mu $ satisfying (2.4) may be complex.

\subhead 3. Regularity in a model case \endsubhead

For a fractional-order $\psi $do $P$ on $\rn$, one can define a homogeneous Dirichlet
problem over an open subset $\Omega \subset \rn$ as follows: For $f$ given
on $\Omega $, find $u$ satisfying
$$
r^+Pu=f \text{ in }\Omega ,\quad u=0 \text{ on }\rn\setminus \Omega ;\tag3.1
$$
this is sometimes called the {\it restricted fractional Dirichlet
problem.} Here $r^+$ indicates
restriction from $\rn$ to $\Omega $; in the cases we study here,  $\Omega =\rnp$.
When $P$ is strongly elliptic, a solvable (perhaps just Fredholm
solvable) problem in low-order $L_2$-Sobolev spaces can be formulated
using variational theory, when $\Omega $ is bounded or $P$ has a
suitable lower boundedness. Then comes the question of {\it regularity}, namely what smoothness can
be obtained for $u$ when $f$ has a given smoothness. Ros-Oton and Serra
gave interesting results for the fractional Laplacian in H\"older spaces in \cite{RS14}, followed
up with generalizations in subsequent papers, and the present author
treated the question in \cite{G15} in $L_p$ Sobolev-type spaces for general $x$-dependent $m$-order pseudodifferential operators $P$ and smooth
domains, satisfying the $\mu $-transmission condition with an
arbitrary complex $\mu $ (followed up in other function spaces in
\cite{G14}).

Recall the notation: $\Cal E'(\rn)$ is the space of distributions with
compact support, $\Cal S(\rn)$ is the Schwartz space of $C^\infty
$-functions $f$ on $\rn$ such that $x^\beta D ^\alpha f$ is bounded
for all $\alpha ,\beta $, $\Cal S'(\rn)$ is its dual space of
temperate distributions. Denote $r^+\Cal S(\rn)=\Cal S(\crnp)$.

The basic strategy of \cite{G15}, when $\Omega =\rnp$, is to compose  $P$
to the left and right with the order-reducing operators $\Xi _-^{\mu -m
}=\OP((\ang{\xi '}-i\xi _n)^{\mu -m })$ and  $\Xi _+^{-\mu 
}=\OP((\ang{\xi '}+i\xi _n)^{-\mu  })$, arriving at a zero-order
operator $Q=\Xi _-^{\mu -m}P\Xi _+^{-\mu }$ satisfying the 0-transmission condition, as studied in the Boutet
de Monvel calculus (cf.\ e.g.\ \cite{B71,G96,S01,G09}), where results are
obtained in standard Sobolev spaces. Strictly speaking, one needs truly pseudodifferential
versions  $\Lambda _-^{\mu -m },\Lambda _+^{-\mu }$ of $\Xi
_-^{\mu -m},\Xi _+^{-\mu }$ to apply the pseudodifferential calculus
properly, but the use of $\Xi _-^{\mu -m},\Xi _+^{-\mu }$ suffices for
some purposes.

Since $\Xi _+^{-\mu }$ is closely linked with
multiplication by $x_n^\mu $ on $\rnp$, one can show that $P$ operates
nicely on functions $u=e^+x_n^\mu v$  with $v\in C^\infty
(\crnp)$; here $e^+$ denotes extension by zero. Moreover, when $u$ solves (3.1), then the
boundary value of $u/x_n^\mu $ for $x_n\to 0+$ exists. To keep notation simple, we now
just recall the smoothest solution space (more general spaces are
recalled in Section 7 below), defined by:
$$
\E_\mu (\crnp)=e^+x_n^\mu C^\infty (\crnp)\text{ when
}\operatorname{Re}\mu >-1; \tag3.2
$$
 and $\E_\mu (\crnp)$
is defined successively as the linear hull of first-order derivatives
of elements of $\E_{\mu +1}(\crnp)$ when $\operatorname{Re}\mu \le -1$
(then distributions supported in the boundary can occur).

It is shown in \cite{G15, Th.\ 4.2, 4.4} that when $P$ is
of order $m$ and satisfies the $\mu $-transmission condition, then
$$
u\in \E_\mu (\crnp)\cap \E'(\rn) \implies r^+ Pu\in  C^\infty
(\crnp)\cap {\bigcap}_s\ol H^s(\rnp),\tag 3.3
$$
and when $P$ moreover is elliptic with factorization index $\mu _0$
($\equiv \mu $ mod 1),
one can conclude the other way: When $u\in \dot H^\sigma
 (\crnp)$ for some $\sigma >\mu _0-\tfrac12$,
$$
 r^+ Pu\in C^\infty
(\crnp)\cap {\bigcap}_s\ol H^s(\rnp)
 \implies u\in \E_{\mu _0}(\crnp).\tag3.4
$$
Our notation for Sobolev spaces is recalled below in the start of Section 4.

For the operator  $L=\OP(\Cal A+i\Cal B)$ considered in Section 2, with a smooth symbol $\Cal
L$
homogeneous of degree $m=2a>0$ for $\xi \ne 0$, $L$ is covered by the
considerations of \cite{G15}, with $\mu =\mu _0=a+\delta $. Namely, we can modify $\Cal L(\xi )$
for $|\xi |\le 1$ to a symbol $\Cal L_1(\xi )$ that is $C^\infty $
also for $|\xi |\le 1$; then the 
remainder
$\Cal R=\OP(\Cal L(\xi )-\Cal L_1(\xi ))$ is smoothing (sends general
distributions into $C^\infty (\rn)$), and has a bounded symbol. Then $r^+L$ sends $\E_\mu
(\crnp)\cap \E'(\rn) $ into $ C^\infty (\crnp)\cap \bigcap_s\ol
H^s(\rnp)$.
The solutions of the
Dirichlet problem with right-hand side in $ C^\infty (\crnp)\cap \bigcap_s\ol
H^s(\rnp)$ are in $\E_\mu (\crnp)$. Thus we have as a special case of \cite{G15}:

\proclaim{Theorem 3.1} Let $0<a<1$, and let $L=\OP(\Cal L(\xi ))$, where $\Cal L(\xi )$
is homogeneous of order $2a$, with even real part
$\Cal A(\xi )$ and odd  imaginary part $i\Cal B(\xi )$,  such that
$\Cal L(\xi )=\Cal A(\xi )+i\Cal B(\xi )$ is $C^\infty $ for
$\xi \ne 0$, and $\Cal A(\xi )>0$ for $\xi \ne 0$. Let
$\mu $ and $\mu ^*$ be defined by $$
\mu =a+\delta ;\quad
\delta =\pi ^{-1}\operatorname{Arctan} \frac {\Cal B(0,1)}{\Cal
A(0,1)}\in \,]-\tfrac12,\tfrac12[\,,\quad \mu ^*=a-\delta  .\tag3.5
$$
Then with $\sigma >\mu  -\tfrac12$, 
$$
\aligned
&u\in \E_\mu (\crnp)\cap \E'(\rn) \implies r^+ Lu\in  C^\infty
(\crnp)\cap {\bigcap}_s\ol H^s(\rnp),\\
 &r^+ Lu\in C^\infty
(\crnp)\cap {\bigcap}_s\ol H^s(\rnp),u\in \dot H^\sigma (\crnp)
 \implies u\in \E_{\mu }(\crnp).
\endaligned\tag3.6
$$
The analogous  result holds for $L^*=\OP(\Cal A(\xi )-i\Cal B(\xi ))$ with $\mu $ replaced by $\mu ^*$.
\endproclaim

There are similar results in spaces based on $L_p$-Sobolev spaces and
other function spaces, see \cite{G15,G14}.

\subhead 4. Integration by parts in the halfspace case \endsubhead

To keep notation simple, we restrict the attention to real $\mu
$; this covers the case of strongly elliptic homogeneous symbols with
a real kernel function. Recall the notation for $L_2$-Sobolev spaces:
$$
\aligned
H^s (\Bbb R^n)
&=\{u\in \Cal S'(\Bbb R^n)\mid \Cal F^{-1}(\langle{\xi }\rangle^s\hat
u)\in L_2 (\Bbb R^n)\},\\
\overline H ^s(\rnp )&=r^+H^s (\Bbb R^n), \text{ the restricted space},\\
\dot H^s (\crnp )&=\{u\in H^s (\Bbb R^n)\mid
\operatorname{supp}u\subset \crnp \}, \text{ the  supported space}.
\endaligned
$$
Here $ \ol H^s(\rnp )$ identifies with the dual space of
$\dot H^{-s}(\crnp )$ for all $s\in\R$ (the duality extending the
$L_2(\rnp)$ scalar product). When $|s|<\frac12$, there is an identification of $\dot
H^s(\crnp )$ with $\ol H^s(\rnp )$ (more precisely $e^+\ol H^s(\rnp )$).
 The trace operator $\gamma _0\colon u\mapsto \lim_{x_n\to 0+}u(x',x_n)$ extends
 to a continuous mapping $\gamma _0\colon \ol H^s(\rnp)\to
 H^{s-\frac12}(\R^{n-1})$ for $s>\frac12$.

Recall that $\Xi _\pm^t=\OP((\ang{\xi '}\pm
i\xi _n)^t)$; here  if homogeneity is important, $\ang{\xi '}$ can be replaced by $[\xi ']$, a
$C^\infty $-function that equals $|\xi '|$ for $|\xi '|\ge 1$ and
ranges in $[\frac12,1]$ when $|\xi '|\le 1$. These operators have the
mapping properties: 
$$
\Xi ^{t }_+\colon \dot H^s(\crnp )\simto
\dot H^{s- t }(\crnp),\quad
r^+\Xi ^{t }_{-}e^+\colon \ol H^s(\rnp )\simto
\ol H^{s- t } (\rnp ),\text{ all }s,t\in\R,\tag4.1
$$
 and for each $t\in{\Bbb R}$, the operators $\Xi ^t _{+}$ and $r^+\Xi
 ^{t }_{-}e^+$, also denoted $\Xi ^t_{-,+}$, identify with each other's adjoints
over $\crnp$. Recall also the simple composition rules (as noted e.g.\
 in \cite{GK93, Th.\ 1.2}):
$$
\Xi ^s_+\Xi ^t_+=\Xi
^{s+t}_+,\quad \Xi ^s_{-,+}\Xi ^t_{-,+}=\Xi
^{s+t}_{-,+}\text{ for $s,t\in\R$. }\tag4.2
$$
Recall from \cite{G15} that for all $\mu $,
$$
\E_\mu (\crnp)\cap \Cal E'(\rn)\subset\Xi _+^{-\mu }e^+(C^\infty
(\crnp)\cap {\bigcap}_s\ol H^s(\crnp)).\tag4.3
$$
A sharpening:  $e^+x_n^{\mu }\Cal S(\crnp)= \Xi
_+^{-\mu }e^+\Cal S(\crnp)$, is shown below in Lemma 6.1.

To show integration-by-parts formulas, we shall use  methods from
\cite{G16} in a simplified 
version. First we have from \cite{G16, Th.\ 3.1}:

\proclaim{Theorem 4.1} Let $0<a<1$, $|\delta |<\frac12$, and $ \mu
=a+\delta $. Let $w\in C^\infty (\crnp)\cap
\bigcap_{s}\ol H^s(\rnp)$, and let
$u'\in \E_\mu (\crnp)\cap
\E'(\rn)$. Denote  $w'=r^+\Xi _+^\mu u'$; correspondingly $u'=\Xi
_+^ {-\mu }e^+w'$ in view of {\rm \cite{G15, Prop.\ 1.7}}. Then
$$
\int _{\rnp}\Xi _-^{ \mu }e^+w\,\partial_n\bar u'\, dx=
(\gamma _0w,\gamma _0w')_{L_2(\R^{n-1})}+(w,\partial_nw')_{L_2(\rnp)}
.\tag 4.4
$$
Here if $\mu \le 0$, the left-hand side is interpreted as the duality
$$
\ang {\Xi ^\mu _{-,+}w,\partial_nu'}_{\ol H^{\frac12+\varepsilon -\mu
}(\rnp),\dot H^{-\frac12-\varepsilon +\mu }(\crnp)},\text{ also equal
to }\int_{\rnp}\Xi
^a_-e^+w\,\overline{\Xi _+^\delta \partial_nu'}\, dx.\tag 4.5
$$

\endproclaim

\demo{Proof}
In view of (4.1),
 $w\in  \ol H^s(\rnp)$ for all
$s$ implies that  
$r^+\Xi _-^{ \mu }e^+w\in   \ol H^{s }(\rnp)$ for all 
$s$.
By (3.3), $w'=r^+\Xi _+^\mu u'\in  C^\infty
(\crnp)\cap \bigcap_s \ol H^s(\rnp)$.

Consider first the {\bf case where $\mu >0$}.
Since $u'\in \E_{\mu }(\crnp)$ with compact support and is continuous
on $\rn$, $\partial_nu'\in
\E_{\mu -1}(\crnp)$ with compact support. Here $x_n^{\mu -1}$ is integrable
over compact sets. Altogether, the integrand
$r^+\Xi
_-^{ \mu }e^+w\, \partial_n\bar u'$ is  the product of $x_n^{\mu -1}$ with a
compactly supported smooth function (on $\crnp$), so the integral is well-defined.

We can also observe that by the identification of $e^+\ol H^t(\rnp)$ and $\dot
H^t(\crnp)$ for $|t|<\frac12$, $e^+w'\in \dot H^{\frac12-\varepsilon
}(\rnp)$ for any $\varepsilon \in \,]0,1[\,$, 
 so
$$
\partial_n u'=\partial_n\Xi _+^{-\mu }e^+w'\in \partial_n\dot
H^{ \mu +\frac12-\varepsilon }(\crnp)\subset\dot
H^{ \mu -\frac12-\varepsilon }(\crnp).
$$

Thus the
integral may be written as the duality
$$
I=\ang{r^+\Xi _-^{ \mu } e^+w,\partial_n 
u'}_{ \ol H^{\frac12- \mu +\varepsilon }(\rnp) ,\dot H^{ \mu -\frac12-\varepsilon
}(\crnp)}.
$$

Note that $r^+\Xi _-^{ \mu }e^+\colon \ol
H^{\frac12+\varepsilon }(\rnp)\to \ol
H^{\frac12- \mu +\varepsilon }(\rnp)$
has the adjoint $\Xi _+^\mu \colon \dot H^{\mu -\frac12-\varepsilon
}(\crnp)\to \dot H^{-\frac12-\varepsilon
}(\crnp)$. This allows to continue the calculation of $I$ as follows:
$$
I=\ang{w,\Xi _+^\mu \partial_n 
u'}_{ \ol H^{\frac12+\varepsilon }(\rnp) ,\dot H^{-\frac12-\varepsilon
}(\crnp)}=\ang{w,\partial_n\Xi _+^\mu u' }=\ang{w,\partial_ne^+w'}.
$$
Here $w'$ itself is a nice function on $\crnp$, but the extension
$e^+w'$ to $\R^n$
has the jump $\gamma _0w'$ at $x_n=0$,
and there holds the formula
$$
\partial_ne^+w'=(\gamma _0 w')(x')\otimes \delta (x_n)+e^+\partial_nw'.\tag4.6
$$
where $\otimes$ indicates a product of distributions with respect to
different variables ($x'$ resp.\ $x_n$). (4.6) is a distribution version
of Green's formula (cf.\ e.g.\ \cite{G96}
(2.2.38)--(2.2.39)); it has been much used in the literature on
boundary value problems (early in the theory in  e.g.\ Seeley \cite{S66}, Boutet de
Monvel \cite{B66}, H\"ormander \cite{H66}). Recall moreover from distribution
theory (cf.\ e.g.\ \cite{G09} p.\ 307) that the ``two-sided'' trace operator $\widetilde\gamma _0\colon
v(x)\mapsto \widetilde \gamma _0 v=v(x',0)$ has the mapping
$\widetilde\gamma _0^*\colon \varphi (x')\mapsto \varphi (x')\otimes
\delta (x_n)$ as adjoint, with continuity properties
$$
\widetilde \gamma _0\colon H^{\frac12+\varepsilon }({\Bbb R}^n)\to H^{\varepsilon }({\Bbb
R}^{n-1}),\quad
\widetilde \gamma _0^*\colon H^{-\varepsilon }({\Bbb R}^{n-1})\to H^{-\frac12-\varepsilon }({\Bbb
R}^{n}), \text{ for }\varepsilon >0. \tag4.7
$$
Here $\widetilde\gamma _0^*\varphi $ is supported in $\{x_n=0\}$,
 hence lies in $\dot H^{-\frac12-\varepsilon }(\crnp)$. We can then write
$$
\partial_ne^+w'=\widetilde \gamma _0^*(\gamma _0
 w')+e^+\partial_nw'\text{ on }\rn.\tag4.8
$$

Since $w\in \ol H^{\frac12+\varepsilon }(\rnp)$, it has an extension
$W\in H^{\frac12+\varepsilon }(\R^n)$ with
$w=r^+W$, and  $\gamma _0 w= \widetilde
\gamma _0W$.
Then
$$
\ang{w,\widetilde \gamma _0^*(\gamma _0  w')}_{ \ol
H^{\frac12+\varepsilon }(\rnp), \dot H^{-\frac12-\varepsilon
}(\crnp)}=\ang{W,\widetilde \gamma _0^*(\gamma _0  w')}_{  H^{\frac12+\varepsilon }({\Bbb R}^n),  H^{-\frac12-\varepsilon }({\Bbb R}^n)};
$$
this is verified e.g.\ by approximating $\widetilde \gamma _0^*(\gamma
_0w')$ in $\dot H^{-\frac12-\varepsilon }$-norm by a sequence of
functions in $C_0^\infty (\rnp)$.  Here we can use (4.7) to write
$$
\multline
\ang{W,\widetilde \gamma _0^*(\gamma
_0  w')}_{H^{\frac12+\varepsilon }({\Bbb R}^n), H^{-\frac12-\varepsilon
}({\Bbb R}^n) }
=\ang{\widetilde \gamma _0W,\gamma
_0  w'}_{H^{\varepsilon }({\Bbb R}^{n-1}), H^{-\varepsilon }({\Bbb
R}^{n-1}) }\\
=\ang{ \gamma _0 w,\gamma
_0  w'}_{H^{\varepsilon }({\Bbb R}^{n-1}), H^{-\varepsilon }({\Bbb
R}^{n-1}) }=(
\gamma _0w,\gamma _0 w')_{L_2({\Bbb R}^{n-1})}.
\endmultline
$$
In the last step we used  that since both $\gamma _0w$ and
$\gamma _0w'$ are in $H^\varepsilon ({\Bbb R}^{n-1})\subset L_2({\Bbb
R}^{n-1}) $, the duality over the boundary is
in fact an $L_2({\Bbb R}^{n-1})$-scalar product. 

 Then finally
$$
\aligned
I&=\ang{w,\partial_ne^+w'}_{ \ol H^{\frac12+\varepsilon }(\rnp) ,\dot
H^{-\frac12-\varepsilon}(\crnp)}\\
&=\ang{w,\widetilde \gamma _0^*(\gamma _0
w')+e^+\partial_nw'}_{ \ol H^{\frac12+\varepsilon }(\rnp) ,\dot
H^{-\frac12-\varepsilon}(\crnp)}\\
&=(\gamma _0w,\gamma _0 w')_{L_2(\R^{n-1})}+\ang{w,e^+\partial_nw'}_{ \ol
H^{\frac12+\varepsilon } ,\dot H^{-\frac12-\varepsilon}}\\
&=(\gamma _0w,\gamma _0 w')_{L_2(\R^{n-1})}+(w,e^+\partial_nw')_{L_2(\rnp)},
\endaligned 
$$
where we used that $w'\in \bigcap_s\ol H^s (\rnp)$. This
shows (4.2) in the case where $\mu >0$.
 
{\bf Case $\mu \le 0$.} The case $\mu \le 0$ can occur only when $a< \frac12$
 and $\delta <0$. As already noted, $r^+\Xi _-^\mu e^+w=\Xi ^\mu
 _{-,+}w$ is in any $\ol H^t(\rnp)$-space, so we can choose $t$ to fit
 with a $\dot H^{-t}(\crnp)$-space for the right-hand factor
 $\partial_nu'$. This factor belongs (cf.\ (4.3)) to
 $$
 \Cal E_{\mu -1}(\crnp)\subset\Xi
 _+^{1-\mu }e^+(C^\infty (\crnp)\cap {\bigcap}_s\ol H^s(\rnp))\subset
 \Xi  _+^{1-\mu }\dot H^{\frac12-\varepsilon }(\crnp)=\dot H^{-\frac12-\varepsilon +\mu }(\crnp),
$$
so we can take $t=\frac12+\varepsilon -\mu $,  showing the first
interpretation in (4.5).

For the second interpretation, note that
$\Xi ^\mu _{-,+}=\Xi ^\delta _{-,+}\Xi ^a_{-,+}$ (cf.\ (4.2)), where
$\Xi ^\delta _{-,+}$ can be transposed to the right-hand side as the
adjoint $\Xi _+^\delta $.
Then the integrand is a product of functions
$\Xi _{-,+}^aw\in C^\infty (\crnp)$ and $\overline{\Xi _+^{\delta
}\partial_nu'}\in \E_{a-1}(\crnp)$, integrable up to the boundary and belonging to $L_1(\crnp)$,  and $I$ may be written
$$
I=\ang{r^+\Xi _-^{ a } e^+w,\Xi _+^{\delta }\partial_n 
u'}_{ \ol H^{\frac12- a +\varepsilon }(\rnp) ,\dot H^{ a -\frac12-\varepsilon
}(\crnp)}=\int_{\rnp}\Xi
^a_-e^+w\,\overline{\Xi _+^\delta \partial_nu'}\, dx.
$$

 We proceed by carrying $\Xi _{-,+}^a$ over to the
right factor, using again (4.2) and the adjoint properties, which gives
$$
I=\ang{w,\Xi _+^\mu \partial_nu'}_{ \ol H^{\frac12+\varepsilon }(\rnp) ,\dot H^{-\frac12-\varepsilon
}(\crnp)}.
$$
From here on, the proof is completed as in the case $\mu >0$.
\qed
\enddemo

An immediate consequence of Theorem 4.1 is the following integration-by-parts result
for a very special operator:

\proclaim{Theorem 4.2}  Let $\mu =a+\delta $, $  \mu ^*=a-\delta $
 with $a,\delta $ as in Theorem {\rm 4.1}, and consider $P=\Xi _{-}^{\mu ^*}\Xi _+^\mu  $. Let
 $u\in \E_\mu (\crnp)\cap \E'(\rn)$ and $u'\in\E_{\mu ^*}(\crnp)\cap
 \E'(\rn)$.
 Then 
$$
\multline
\int _{\rnp}Pu\,\partial_n\bar u'\, dx+\int
_{\rnp}\partial_nu\,\overline
{P^* u'}\, dx\\
=\Gamma (\mu +1){\Gamma (\mu ^*+1)}
\int_{{\Bbb
R}^{n-1}}\gamma _0(u/x_n^{\mu  })\,\gamma _0(\bar u'/{x_n^{\mu ^*}})\, dx'.
\endmultline
\tag 4.9
$$
\endproclaim

\demo{Proof}
We  apply Theorem 4.1 to the integrals in the left-hand side
of (4.9).
When $u\in \E_\mu (\crnp)\cap \E'(\rn)$, then  $w=r^+\Xi _+^{\mu }u\in \bigcap_s\ol
H^s (\rnp) $. Hence $r^+Pu=r^+\Xi _-^{\mu ^*}e^+w\in \bigcap_s\ol
H^s (\rnp) $, and an application of Theorem 4.1 with $ \mu $ replaced
by $\mu ^*$ gives:
$$
\int _{\rnp}Pu\,\partial_n\bar u'\, dx=(\gamma _0w,\gamma _0w')_{L_2(\R^{n-1})}+(w,\partial_nw')_{L_2(\rnp)},
$$
where $w'=r^+\Xi _+^{ \mu ^*}u'$.

Using that $P^*=\Xi _-^{ \mu }\Xi _+^{ \mu ^*}$, we can apply
the analogous argument to show that the conjugate of $\int
_{\rnp}\partial_nu\,\overline{P^*u'}\,dx$ satisfies
$$
\int _{\rnp}P^*u'\,\partial_n\bar u\, dx=(\gamma _0w',\gamma _0w)_{L_2(\R^{n-1})}+(w',\partial_nw)_{L_2(\rnp)},
$$
with the same definitions of $w'$ and $w$.
It follows by addition that
$$
\multline
\int _{\rnp}Pu\,\partial_n\bar u'\, dx+\int
_{\rnp}\partial_nu\,
\overline {P^* u'}\, dx\\
=2(\gamma _0w,\gamma _0w')_{L_2(\R^{n-1})}+\ang{w,\partial_n 
w'}_{L_2(\rnp)}+\ang{\partial_nw, 
w'}_{L_2(\rnp)}=(\gamma _0w,\gamma _0w')_{L_2(\R^{n-1})};
\endmultline
$$
in the last step we used that $\int_{\rnp}(w\partial_n\bar w'+\partial_nw \bar
w')\,dx=-\int_{{\Bbb R}^{n-1}}\gamma _0w\gamma _0\bar w'\, dx'$.
Insertion of $\gamma _0w=\gamma _0\Xi _+^{\mu }u=\Gamma(1+\mu )\gamma
_0(u/x_n^{\mu })$ (as shown in \cite{G15, Th.\ 5.1} and reproved below
in (6.9)
for $\mu -1$) and  $\gamma _0w'=\gamma _0\Xi _+^{\mu ^* }u'
 =\Gamma(1+\mu ^*)\gamma
_0(u'/x_n^{\mu ^* })$ gives (4.9).
\qed
\enddemo


We now turn to the operator 
 $L=\OP(\Cal A(\xi )+i\Cal B(\xi ))$ 
 considered in 
 Theorem 3.1. 
Define
$$
\aligned
Q=\Xi _-^{-\mu ^*}L\Xi _+^{-\mu }=\OP(q(\xi )),&\quad q(\xi )=(\ang{\xi
'}-i\xi _n)^{-\mu ^*}\Cal L(\xi )(\ang{\xi '}+i\xi _n)^{-\mu },\\
\text{so that }L=\Xi _-^{\mu ^*}Q\Xi _+^{\mu },&\quad \Cal L(\xi )=(\ang{\xi
'}-i\xi _n)^{\mu ^*}q(\xi )(\ang{\xi '}+i\xi _n)^{\mu }.
\endaligned
\tag4.10
$$
Since $\Cal L$ is $C^\infty $ only for $\xi \ne 0$, and just H\"older
continuous of order $2a$ at 0, we shall write it as the sum of a
smooth symbol and a symbol with small support:
$$
\Cal L(\xi )=\Cal L_\varphi  (\xi )+\varphi  (\xi )\Cal
L(\xi ),
$$
where $\varphi (\xi )$ is a  function in $C_0^\infty (\rn,[0,1])$
supported for $|\xi |\le 1 $ and equal to 1 for $|\xi |\le
\frac12$; so that $\Cal L_\varphi  (\xi )=(1-\varphi
 (\xi ))\Cal L(\xi )$ is in $S^{2a}(\rn\times\rn)$.
Then similarly,
$$
q(\xi )=q_\varphi  (\xi )+\varphi  (\xi )q(\xi ),
$$
where $q_\varphi (\xi )$ is a smooth symbol of order 0. Its principal
part $q^0=(|\xi '|-i\xi _n)^{-\mu ^*}\Cal L(\xi )(|\xi '|+i\xi _n)^{-\mu }$ is homogeneous of
degree 0 in $\xi $ and satisfies the
$0$-transmission condition:
$$
\aligned
q^0(0,1)&=(-i)^{-a+\delta }\Cal L(0,1)i^{-a-\delta }=i^{-2\delta
}\Cal L(0,1)=e^{-i\pi \delta }\Cal L(0,1),\\
q^0(0,-1)&=(+i)^{-a+\delta }\Cal L(0,-1)(-i)^{-a-\delta }=i^{2\delta
}e^{i\pi (2a-2(a+\delta )}\Cal L(0,1)\\
&=e^{-i\pi \delta }\Cal L(0,1)=q^0(0,1).
\endaligned
\tag 4.11
$$
Note that $$
s_0\equiv q^0(0,1)=e^{-i\pi \delta }\Cal
L(0,1)=|\Cal L(0,1)|.
\tag4.12
$$
It is checked by use of a Taylor expansion of $\ang{\xi '}=|\xi '|(1+|\xi '|^{-2})^\frac12$ that also the lower-order terms in $q_\varphi  $ fulfill
the rules for the 0-transmission condition (cf.\ \cite{G15}).

 \proclaim{Theorem 4.3} Consider  $L=\OP(\Cal L(\xi ))$, $a$, $\delta
 $, $\mu $, $\mu ^*$  as described
 in Theorem {\rm 3.1}. 
 For $u\in \E_\mu (\rnp)\cap \E'(\rn)$, $u'\in \E_{\mu ^*} (\rnp)\cap \E'(\rn)$,
 there holds 
$$
\aligned 
\int_{\rnp} Lu\,\partial_n\bar
u'\,dx&+\int_{\rnp}\partial_nu\,\overline{  L^*u'}\,dx\\
&=\Gamma (\mu +1){\Gamma(\mu ^*+1)}\int_{{\Bbb
R}^{n-1}}s_0\gamma _0(u/x_n^{\mu })\,{\gamma _0(\bar u'/x_n^{\mu ^*})}\,
dx',
\endaligned
\tag4.13
$$
where $s_0=|\Cal L(0,1)|=(\Cal A(0,1)^2+\Cal B(0,1)^2)^\frac12$.

The integrals over $\rnp$ are interpreted as in Theorem {\rm 4.1} if $\mu $ or
$\mu ^*\le 0$.
\endproclaim

\demo{Proof}
We first consider the contribution from $\varphi \Cal
L$. This can be handled in the same way as in the treatment of the
smoothing term $\Cal S$ in the proof of \cite{G16, Cor.\ 3.5}, showing
that  $\varphi \Cal L$ contributes with
 zero.

Now consider the contribution from $\Cal L_\varphi  $. To simplify
the notation we omit the subscript $\varphi $ from now on.
 It is accounted for in \cite{G16} that $q_\varphi $ (from here on
 simplified to $q$) 
 lies in $S^0(\R^{n-1}, \Cal H_0)$, hence has a decomposition
$$
q(\xi ',\xi _n)=q(0,1)+f_+(\xi )+f_-(\xi ),\quad q(0,1)=s_0\tag4.14
$$
(cf.\ (4.12)), where
$$
f_+(\xi ',\xi _n)=h^+q\in S^0(\R^{n-1}, \Cal H^+),\quad f_-(\xi ',\xi _n)=h^-_{-1}q\in S^0(\R^{n-1}, \Cal H^-_{-1}).
$$
The teminology with $h^\pm$ and $\Cal H_d$, $\Cal H^\pm$, etc., originating from
a decomposition principle ascribed to Wiener and Hopf, and set up in
this way in Boutet de Monvel \cite{B71},  is summarized
in \cite{G16} and in \cite{G18}, and described at length e.g.\  in
\cite{G96}, \cite{G09} and Schrohe \cite{S01}. Let
us here just recall that the functions $f(\xi _n)\in \Cal H_0$ are the
Fourier transforms of $\tilde f(x_n)\in e^+\Cal S(\crp)\oplus e^-\Cal
S(\crm)\oplus \C\delta =\Cal F^{-1}(\chp\oplus \chm_{-1}\oplus \C) $, where
$h^+$ corresponds to the projection onto $e^+\Cal S(\crp)$, and
$h^-_{-1}$ corresponds to the projection onto $e^-\Cal S(\crm)$. Here
$\Cal S(\crpm)=r^\pm\Cal S(\R)$, the Schwartz space. (More details in
\cite{G18, (A.2)ff.}.)

In \cite{G16} it is moreover shown that $q$ has a {\it
factorization} into two such factors, by making use of the strong
ellipticity that allows passing to $\log q$ and back,
but we do not need the factorization nor the ellipticity here, only
the sum decomposition (4.14). This simplifies the proof and is inspired from \cite{G18}.

Denote
$$
\OP (f_+)=F_+,\quad \OP(f_-)=F_-,\text{ so }F_-^*=\OP(\overline f_-),\tag4.15
$$
where $\overline f_-\in S^0(\R^{n-1}, \Cal H^+)$. Then
$$
L=\Xi _-^{-\mu ^*}Q\Xi _+^{\mu },\quad Q=\OP(q)=s_0+F_++F_-.\tag4.16
$$
The contribution from the first term in $Q$ is $s_0\Xi _-^{-\mu ^*}\Xi
_+^{\mu }$ which we have already treated in Theorem 4.2; it gives the
right-hand side in (4.13).

It remains to treat the terms with $F_+$ and $F_-$ and to show that
they together contribute with 0. Here we note that as in \cite{G15,
Th.\ 4.2}, using that $u$ is supported in $\crnp$,
$$
r^+ \Xi _-^{\mu ^*}F_+\Xi _+^{\mu }u= r^+\Xi _-^{\mu ^*}e^+(F_+\Xi
_+^{\mu })u,\quad
r^+ \Xi _-^{\mu ^*}F_-\Xi _+^{\mu }u= r^+(\Xi _-^{\mu ^*}F_-)e^+\Xi _+^{\mu }u.\tag4.17
$$
(In the
various calculations, the extension by 0 on $\rnm$ is sometimes tacitly
understood.)

For the given $u,u'$, denote
$$
\aligned
w&=r^+\Xi _+^\mu u,\quad w_1=r^+F_+w=r^+F_+\Xi _+^\mu u,\\
 w'&=r^+\Xi _+^{\mu ^*}u',\quad w'_2=r^+F^*_-w'=r^+F_-^*\Xi _+^{\mu
 ^*}u';
\endaligned
\tag4.18 
$$
they all lie in $C^\infty (\crnp)\cap \bigcap_s\ol H^s(\rnp)$. 

For the first term in (4.13), we proceed as follows: By Theorem 4.1
(applied with $\mu $ replaced by $\mu ^*$), we have for the
contribution from $F_+$, in
view of (4.17),
$$
\multline
\int_{\rnp}\Xi _-^{\mu ^*}F_+\Xi ^\mu _+u\partial_n\bar u'\,dx=
\int_{\rnp}\Xi _-^{\mu ^*}e^+w_1\partial_n\bar u'\,dx\\
=(\gamma _0w_1,\gamma
_0w')_{L_2(\R^{n-1})}+(w_1,\partial_nw')_{L_2(\rnp)}.
\endmultline
\tag4.19
$$
We can likewise use the argumentation of Theorem 4.1 to show for the
contribution from $F_-$, after a transposition:
$$
\multline
\int_{\rnp}\Xi _-^{\mu ^*}F_-\Xi ^\mu _+u\partial_n\bar u'\,dx=
\int_{\rnp}\Xi ^\mu _+u\overline{(\Xi _-^{\mu ^*}F_-)^*\partial_n u'}\,dx\\
=\int_{\rnp}\Xi ^\mu _+u\overline{F_-^*\Xi _+^{\mu ^*}\partial_n u'}\,dx
=(\gamma _0w,\gamma
_0w'_2)_{L_2(\R^{n-1})}+(w,\partial_nw'_2)_{L_2(\rnp)}.
\endmultline
$$

Now there is a special observation that allows removing the boundary
terms appearing here, similarly as at \cite{G16,(3.30)}:
With $w_1$  defined in (4.18) we have, using that the boundary value of a function
supported in $\crnp$ can be described by the convention for trace
operators in the Boutet de Monvel calculus (\cite{G18, (A.15), (A.1)}):
$$
\gamma _0w_1=\gamma _0(F_+w)=(2\pi )^{-n}\int_{{\Bbb
R}^{n-1}}e^{ix'\cdot\xi '}\int_{\R}f_+(\xi ',\xi _n)\F(e^+w)\,d\xi
'd\xi _n.
$$
This equals 0 for the following reason:
Both $f_+$ and $\F(e^+w)$ are in $\chp$ as
functions of $\xi _n$, in particular $O(\ang{\xi _n}^{-1})$, whereby the integrand is $O(\ang{\xi _n}^{-2})$ and extends holomorphically into the
lower imaginary halfplane $\C_-$; then the integral over $\R$ can
be transformed to a closed contour in $\C_-$ and therefore vanishes.

Similarly, $\gamma _0w'_2=0$. Thus only the integrals over $\rnp$
remain, and we find altogether:
$$
\int_{\rnp}\Xi _-^{\mu ^*}(F_++F_-)\Xi ^\mu _+u\,\partial_n\bar
u'\,dx=(w_1,\partial_nw')_{\rnp}+(w,\partial_nw'_2)_{\rnp}.\tag4.20
$$

 The analogous arguments apply to the contribution from $F_++F_-$ in
 the second term in (4.13) (after
 conjugation), with $\mu $ and $\mu ^*$ interchanged, and lead to:
$$
\int_{\rnp} \partial_nu\,\overline{(\Xi _-^{\mu ^*}(F_++F_-)\Xi ^\mu _+)^*
u'}\,dx=(\partial_nw,w'_2)_{\rnp}+(\partial_nw_1,w')_{\rnp}.\tag4.21
$$

Addition of the right-hand sides of (4.20) and (4.21) gives
$$
\multline
(w_1,\partial_nw')_{\rnp}+(w,\partial_nw'_2)_{\rnp}
+(\partial_nw,w'_2)_{\rnp}+(\partial_nw_1,w')_{\rnp}\\
=\int_{\rnp}(\partial_n(w\bar w'_2)+\partial_n(w_1\bar w'))\, dx
=-\int_{\R^{n-1}}(\gamma _0w\gamma _0\bar w'_2+\gamma _0w_1\gamma
_0\bar w')\, dx'=0,
\endmultline
$$
where we have again used the vanishing of $\gamma _0w_1$ and $\gamma
_0w'_2$.

Adding all the contributions from $F_+$, $F_-$ and $s_0$, we arrive at (4.13).\qed
\enddemo

Theorem 4.3 is proved by real methods in \cite{DRSV20, Prop.\ 5.1}, for operators $L$
that in addition satify (2.9), and are  infinitesimal
generators of  $\alpha $-stable $n$-dimensional L\'evy processes. It
can be remarked that in the above method, we only work with the
{\it symbol} $\Cal L(\xi )$ of $L$ and not with the corresponding {\it
distribution kernel},
so that we avoid the issue of singularities of the kernel.

The results can undoubtedly be extended to $x$-dependent operators by
efforts as in \cite{G16}, as long as $\mu $ is constant in $x'$. For
variable $\mu (x') $, other tools are needed. There is a  method in \cite{DRSV20}
showing how one extends the integration-by-parts formula from the halfspace case to general domains for the
$x$-independent operators considered there.

\subhead 5. The nonhomogeneous Dirichlet problem and a halfways Green's formula \endsubhead

Along with the homogeneous Dirichlet problem (3.1), one can consider a
nonhomogeneous Dirichlet problem if  the scope expanded to allow
so-called "large solutions", behaving like $d^{\mu -1}$ near the
boundary (with $d(x)=\operatorname{dist}(x,\partial\Omega )$); such solutions blow up at the boundary when $\mu <1$.

For elliptic operators $P$ satisfying the $\mu
$-transmission condition, we
can then pose the nonhomogeneous problem 
$$
r^+Pu=f \text{ in }\Omega ,\quad \gamma _0(u/d^{\mu -1} )=\psi 
\text{ on }\partial\Omega ,\quad u=0 \text{ on }\rn\setminus \Omega ;\tag5.1
$$
assuming $\mu >0$. It has good solvability properties when $u$ is sought in $\Cal
E_{\mu -1} (\comega)$ and related Sobolev-type spaces, cf.\
\cite{G15, Sect.\ 5 and Th.\ 6.1}. (We now take  $\mu > 0$, since $\E_{\mu
-1}(\comega)$ is defined via (3.2) then; for lower $\mu $, the
interpretation of the boundary condition may be more complicated.)
The homogeneous Dirichlet problem fits in here as the problem with
$\psi =0$ in (5.1);  for $f\in C^\infty
(\comega)$ the solutions  with
vanishing $\gamma _0(u/d^{\mu -1})$ are in $\Cal E_\mu
(\comega)$.

The interest of problem (5.1) for the fractional
Laplacian $(-\Delta )^a$ (where $\mu =a$) was also pointed out in
Abatangelo \cite{A15} (independently of
\cite{G15});  the boundary condition there is given
in a less explicit way, except when $\Omega $ is a ball. 

One can now set up a ``halfways Green's formula'' describing the
difference of 
integrals of $Pu$ times $v$ and $u$  times $P^*v$ as a boundary
integral, when $u$ is a large solution for $P$ (say in $\Cal E_{\mu
-1}(\comega)$) and $v$ is an ordinary solution for $P^*$. Such a
formula was shown in \cite{G18} when $P$ is even (so that  $\mu =a$), and a related
formula was shown in \cite{A15, formula 1.2(9)} for $(-\Delta )^a$. 

It was observed in \cite{G18} for even operators
that the halfways Green's formula is
essentially equivalent with the integration-by-parts formula (in Cor.\
4.5 there, the conclusion in one direction is explained; in fact the same
ingredients allow the other conclusion, as also done below).

We shall generalize this to the present operators, and can thereby
rather easily deduce a halfways Green's formula in the
case of $L$ on a halfspace, by manipulations with the result of
Theorem 4.3. The word ``halfways'' refers to taking one function $u$
as a solution of the nonhomogeneous Dirichlet problem (for $L$) and
the other function $v$  as a solution of the homogeneous Dirichlet problem
(for $L^*$). (When both functions are solutions of nonhomogeneous Dirichlet
problems, one should get a full Green's formula with more boundary
terms; this is carried out in Section 6 below.)

\proclaim{Theorem 5.1} Let  $L=\OP(\Cal L(\xi ))$, $a$, $\delta $,
 $\mu $, $\mu
 ^*$ and $s_0$ be as described in Theorems {\rm 3.1} and {\rm 4.3},
 and assume moreover that $\mu , \mu ^* >0$.
 
 For $u\in \E_{\mu -1}(\crnp)\cap \E'(\rn)$ and  $v\in \E_{\mu ^*} (\crnp)\cap \E'(\rn)$,
 there holds 
$$
\int_{\rnp} Lu\,\bar v\,dx-\int_{\rnp}u\,\overline{  L^*v}\,dx=-\Gamma (\mu ){\Gamma(\mu ^*+1)}\int_{{\Bbb
R}^{n-1}}s_0\gamma _0(u/x_n^{\mu -1 })\,{\gamma _0(\bar v/x_n^{\mu ^*})}\,
dx'.
\tag5.2
$$

\endproclaim

As shown for the fractional Laplacian in \cite{A15}, this type of formula can be used to derive a representation of a solution of
(5.1) in terms of the data.
 
Note that since $\mu =a+\delta $ and $\mu ^*=a-\delta $ with $|\delta
|<\frac12$, the extra requirement on positivity of $\mu ,\mu ^*$ is always
satisfied when $a\ge \frac12$.

As a preparation for the proof, we recall an elementary observation:

\proclaim{Lemma 5.2} Let $\mu >0$. Then any function $u\in \E_{\mu
-1}(\crnp)\cap \E'(\rn)$  can be written as $$
u=\partial_nU+u_1, \text{ where }U,u_1\in \E_{\mu}(\crnp)\cap
\E'(\rn).\tag 5.3
$$
\endproclaim 

There is a detailed proof of this on page 494 in \cite{G15}. We also
note an easy identity:

\proclaim{Theorem 5.3}  Consider  $L=\OP(\Cal L(\xi ))$, $a$, $\delta
 $, $\mu $, $\mu ^*$  as described
 in Theorem {\rm 3.1}. 
 For  $w,w'\in \dot H^a(\crnp)$,
there holds
$$
\ang{r^+Lw, w'}_{\ol H^{-a}(\rnp),\dot H^a(\crnp)}- \ang{ w,r^+L^*w'}_{\dot H^{a}(\crnp),\ol H^{-a}(\rnp)}=0. \tag5.4
$$

In particular,
$$
\int_{\rnp} Lw\,\bar w'\,dx-\int_{\rnp}w\,\overline{
L^*w'}\,dx=0 \tag5.5
$$
holds when
$w\in \E_{\mu}(\crnp)\cap
\E'(\rn) $, $w'\in \E_{\mu ^*}(\crnp)\cap \E'(\rn) $.
\endproclaim

\demo{Proof} For $w,w'\in H^a(\rn)$ we have
$$
\ang {Lw,w'}_{H^{-a}(\rn), H^a(\rn)}-\ang {w,L^*w'}_{H^{a}(\rn), H^{-a}(\rn)}=0,
$$
since $L$ and its adjoint $L^*$ are bounded from  $H^{a}(\rn)$ to $
H^{-a}(\rn)$ (having $x$-independent symbols that are $O(|\xi
|^{2a})$). When moreover $w,w'\in \dot H^a(\crnp)$,
$$
\aligned
\ang {Lw,w'}_{H^{-a}(\rn), H^a(\rn)}&=\ang {r^+Lw,w'}_{\ol H^{-a}(\rnp),
\dot H^a(\crnp)},\\
\ang {w,L^*w'}_{H^{a}(\rn), H^{-a}(\rn)}&=\ang {w,r^+L^*w'}_{\dot
H^{a}(\crnp),\ol H^{-a}(\rnp)};
\endaligned
$$ 
e.g.\ for the first expression, we can approximate $w'$ in $\dot
H^{a}(\crnp)$ by functions $\varphi \in C_0^\infty (\rnp)$, where it
clearly holds. Thus (5.4) holds for $w,w'\in \dot H^a(\crnp)$.

According to \cite{G15, Prop.\ 4.1}, $ \E_{\mu}(\crnp)\cap
\E'(\rn)\subset H^{\mu (s)}(\crnp)$ for any $s>\mu -\frac12$;  here
$H^{\mu (s)}(\crnp)=\dot H^s(\rnp)$ if $s\in \,]\mu -\frac12,\mu
+\frac12[\,$ (cf.\ \cite{G15, Def.\ 1.4 ff.\ or Th.\ 5.4}). Since $\mu
=a+\delta $ with $|\delta |<\frac12$, the value $s=a$ satisfies
$$
\mu -\tfrac12=a+\delta -\tfrac12<a<a+\delta +\tfrac12=\mu +\tfrac12,
$$
so $H^{\mu (a)}(\crnp)=\dot H^a(\rnp)$. Thus $ \E_{\mu}(\crnp)\cap
\E'(\rn)\subset \dot H^a(\rnp) $. Similarly, since $\mu ^*=a-\delta $
with   $|\delta |<\frac12$,  $ \E_{\mu^*}(\crnp)\cap
\E'(\rn)\subset \dot H^a(\rnp) $. Then formula (5.4) is valid for such
functions, where it may be written in terms of integrals (5.5).\qed
\enddemo

\demo{Proof of Theorem {\rm 5.1}}
Use Lemma 5.2 to write the given $u$ as in (5.3),
$u=\partial_nU+u_1$.

The contribution from $u_1$ is dealt with by Theorem 5.3, showing that
$$
\int_{\rnp} Lu_1\,\bar v\,dx-\int_{\rnp}u_1\,\overline{  L^*v}\,dx=0. \tag5.6
$$
For the contribution from $\partial_nU$, we note that, writing
$U=x_n^\mu w$ for $x_n>0$, $w\in C^\infty (\crnp)$,
$$
\partial_nU=\partial_n(x_n^\mu w)=\mu x_n^{\mu
-1}w+x_n^\mu \partial_nw 
\text{ for }x_n>0,
$$
so the weighted boundary value for $x_n\to 0+$ satisfies
$$
\gamma _0(\partial_nU/x_n^{\mu -1})=\mu \gamma _0w=\mu \gamma _0(U/x_n^{\mu }).\tag5.7
$$
Now by Theorem 4.3 applied to $U$ and $v$,
$$
\aligned
\Gamma (\mu +1)&{\Gamma(\mu ^*+1)}\int_{{\Bbb
R}^{n-1}}s_0\gamma _0(U/x_n^{\mu })\,{\gamma _0(\bar v/x_n^{\mu
^*})}\,dx'
=\int_{\rnp} LU\,\partial_n\bar
v\,dx+\int_{\rnp}\partial_nU\,\overline{  L^*v}\,dx\\
&=-\int_{\rnp}\partial_n( LU)\,\bar
v\,dx+\int_{\rnp}\partial_nU\,\overline{  L^*v}\,dx
=-\int_{\rnp}L(\partial_n U)\,\bar
v\,dx+\int_{\rnp}\partial_nU\,\overline{  L^*v}\,dx;
\endaligned
$$
 here we used in the
 integration by parts that $\gamma _0v=0$ (since $\mu ^*>0$). 
 In view
 of (5.7), this shows that 
$$
\int_{\rnp}L(\partial_n U)\,\bar
v\,dx-\int_{\rnp}\partial_nU\,\overline{  L^*v}\,dx=-\Gamma (\mu ){\Gamma(\mu ^*+1)}\int_{{\Bbb
R}^{n-1}}s_0\gamma _0(\partial_nU/x_n^{\mu -1})\,{\gamma _0(\bar v/x_n^{\mu
^*})}\,dx'.
$$
Adding (5.6) to this, we find the desired formula (5.2). \qed 
\enddemo

One can conversely derive Theorem 4.3 from Theorem 5.1.

\subhead 6. The full Green's formula \endsubhead

The formulas in Theorem 4.1 and 5.3 were shown by very basic tools in
 Fourier theory and distribution theory, and are surprising by containing only local terms
in the integrals over the boundary, both in the sense that the trace
operators $\gamma _0^{\mu -1}$ and $\gamma _0^\mu $ are local, and
that they enter with the local coefficient $s_0$. For completeness, we shall now
show a full Green's formula (of which they are corollaries). Here
we need to draw on some further tools from the Boutet de Monvel
calculus, and the resulting integrals over the boundary contain local trace
operators but also a nonlocal coefficient in the form of a first-order
pseudodifferential operator over the boundary. The proof 
follows the scheme worked out in 
\cite{G18} for the case $\mu =a$, of course with simplifications due to considering only a model
problem.

Recall the notation from \cite{G15} and \cite{G18, Sect.\ 3}
\footnote{In \cite{G18, Sect.\ 3}, the indication $ e^+\Cal S(\crnp)$ 
 should be replaced by $ e^+x_n^{a -1}\Cal S(\crnp)$ in the
 occurrences on page 756, 757 and 758.}:
When
$u\in  e^+x_n^{\mu -1}\Cal S(\crnp)$ with $\mu >0$, we have an
expansion at the boundary defined by Taylor expanding $w=u/x_n^{\mu
-1}$ for $x_n\ge 0$ and
renormalizing the coefficients:
$$
u(x)=u_0(x')I^{\mu -1}(x_n)+u_1(x')I^{\mu }(x_n)+\dots+ u_k(x')I^{\mu
-1+k}(x_n)+O(x_n^{\mu +k}),\tag6.1
$$
where  $I^\mu (x_n)$ is defined as in \cite{G15}:
$$
I^\mu (x_n)= H(x_n)x_n^\mu /\Gamma (\mu +1)\text{ when
}\operatorname{Re}\mu >-1, \text{ here } \partial_{x_n}I^{\mu +1}=I^{\mu },\tag6.2
$$
$H$ being the Heaviside function $1_{\{x_n\ge 0\}}$. The
Gamma factor serves to normalize $I^\mu
$  so that $I^{\mu }=\partial_{x_n}I^{\mu +1}$; this formula is also
used to define the distribution $I^\mu $ for lower
$\operatorname{Re}\mu $. 

 The coefficients $u_k(x')\in \Cal S(\R^{n-1})$ in (6.1)
are denoted $\gamma _k^{\mu -1}u$; here 
$$
\aligned
&\gamma _k^{\mu -1}u=\Gamma (\mu +k)\gamma _k(u/x_n^{\mu -1}),\text{
in particular}\\
&\gamma _0^{\mu -1}u=u_0=\Gamma (\mu )\gamma
_0(u/x_n^{\mu -1}),\quad
 \gamma _1^{\mu -1}u=u_1=\Gamma (\mu +1)\gamma
 _1(u/x_n^{\mu -1}),
 \endaligned\tag6.3
$$
where $\gamma
_kw=\lim_{x_n\to 0+}\partial_n^kw$.
The cases $k=0$ and 1 are
viewed as the {\it Dirichlet 
resp.\ Neumann
traces} of $u\in \Cal E_{\mu -1}(\crnp)$. 
It is useful to involve also another expansion obtained by Taylor
expanding $x_n^{-\mu }e^{\ang{\xi '}x_n}\F_{x'\to\xi '}u$; here the partially Fourier transformed 
terms have a factor $e^{-\sigma x_n}$, $\sigma =\ang{\xi '}$:
$$
\aligned
\Cal F_{x'\to\xi '}u&\equiv \acute u(\xi ',x_n)
=\hat \varphi _0(\xi ')I^{\mu -1}(x_n) e^{-\sigma x_n}+\acute u'(\xi ',x_n)\\
&=\hat
\varphi _0(\xi ')I^{\mu -1}(x_n) e^{-\sigma x_n}+\hat \varphi _1(\xi ')I^{\mu }(x_n) e^{-\sigma
x_n}+\acute u''(\xi ',x_n)\\
&\sim\sum_{k\ge 0}\hat \varphi _k(\xi ')I^{\mu -1+k}(x_n)e^{-\sigma
x_n},
\endaligned\tag6.4
$$
where the $\varphi _k$ are in $\Cal S(\R^{n-1})$. Furthermore,
$$
\Cal Fu\sim \sum_{k\ge 0}\hat \varphi _k(\xi ')(\sigma +i\xi _n)^{-\mu -k},\tag6.5
$$
in view of the formula 
$$
\Cal F_{x_n\to \xi _n}[I^\mu(x_n) e^{-\sigma x_n}]=(\sigma +i\xi _n)^{-\mu
-1}.\tag6.6
$$
(The two expansions (6.1) and (6.5)
are examined in \cite{G15, Sect.\ 5} as a tool for the study of
nonhomogeneous boundary problems.)
By the expansion (6.4), we can write
$$
\aligned
u&=U_0+u'=U_0+U_1+u''\sim \sum_{k\ge 0}U_k,\text{ with }\\
U_k&=\Cal F^{-1}_{\xi '\to x'}[\hat \varphi _k(\xi ')I^{\mu -1+k}(x_n)e^{-\sigma
x_n}]=\Cal F^{-1}[\hat\varphi _k(\xi ')(\sigma +i\xi _n)^{-\mu -k}].
\endaligned
\tag6.7
$$
Since $\F^{-1}_{\xi '\to x'}[\hat \varphi _k(\xi ')H(x_n)e^{-\ang{\xi '}
x_n}]\in e^+\Cal S(\crnp)$, one has that $U_0\in e^+x_n^{\mu -1}\Cal S(\crnp)$, $U_1\in
e^+x_n^\mu \Cal S(\crnp) $, and $U_k\in
e^+x_n^{\mu +k-1}\Cal S(\crnp) $ in general.

There is a one-to-one correspondence between the coefficient sets
$\{u_0,\dots,u_k\}$ and $\{\varphi _0,\dots,\varphi _k\}$ for any $k$, that
follows by comparing the Taylor expansions of $\acute w=$\linebreak$\F_{x'\to\xi
'}(u/x_n^{\mu -1})$  and $\acute w_e=\F_{x'\to\xi
'}(e^{\sigma x_n}u/x_n^{\mu -1})$. We just need the
transition formula for the first two coefficients, namely, when the
$I^\mu $-factors are taken into account:
$$\aligned
\hat \varphi _0&=\Gamma (\mu )\gamma _0\acute w_e=\Gamma (\mu )\gamma
_0\acute w= \hat u_0,\\
\hat \varphi _1&=\Gamma (\mu +1)\gamma _0(\partial_n\acute w_e)=\Gamma (\mu +1)\gamma
_0(\sigma \acute w_e+e^{\sigma x_n}\partial_n\acute w)\\
&= \Gamma (\mu +1)\Gamma (\mu )^{-1}\sigma \hat u_0+\Gamma (\mu
+1)\gamma _1\acute w=\mu \sigma \hat u_0+\hat u_1;
\endaligned$$
hence with  $\ang{D'}=\OP(\ang{\xi '})$,
$$
\varphi _0=u_0,\quad \varphi _1=u_1+\mu \ang{D'}u_0.\tag6.8
$$

Note that since $\gamma _0^{\mu -1}U_0=\varphi _0=u_0=\gamma _0^{\mu
-1}u$, the first trace of $u'=u-U_0$ is zero, so  $u'\in e^+x_n^\mu \Cal S(\crnp) $. Similarly, having the two
first traces equal to zero,  $u''=u-U_0-U_1$ is in $ e^+x_n^{\mu
+1}\Cal S(\crnp)$. In general, $u^{(k)}=u-(U_0+\dots+U_{k-1})$ is in $ e^+x_n^{\mu
+k-1}\Cal S(\crnp)$.

We also observe:
$$
\gamma _0^{\mu -1}u=\gamma _0\Xi _+^{\mu -1}u.                               \tag6.9
$$
This follows since
$$
\aligned
\gamma _0\Xi _+^{\mu -1}u&=\gamma _0\Xi _+^{\mu -1}(U_0+\dots
+U_k+u^{(k)}),\text{ where}\\
\gamma _0\Xi _+^{\mu -1}U_j&=\gamma _0\Cal F^{-1}((\sigma +i\xi
_n)^{\mu -1}\hat\varphi _j(\sigma +i\xi _n)^{-\mu -j})=\gamma _0\Cal
F^{-1}(\hat\varphi _j(\sigma +i\xi _n)^{-1 -j})\\
&=\varphi _0\text{ if
}j=0,\;0\text{ if }j>0,\\
\Xi _+^{\mu -1}u^{(k)}&\in \Xi _+^{\mu -1}\dot H^{k-2}(\crnp)=\dot
H^{k-1-\mu }(\crnp)\text{ for large }k; \text{ then }\gamma _0\Xi _+^{\mu -1}u^{(k)}=0.
\endaligned
$$

Note also: 
$$
\aligned
\Xi _+^\mu U_k&=\Cal F^{-1}[(\sigma +i\xi _n)^\mu \hat \varphi _k(\xi ')(\sigma +i\xi
_n)^{-\mu -k}]\\
&=\Cal F^{-1}_{\xi '\to x'}[\hat \varphi _k(\xi
')I^{k-1}e^{-\sigma x_n}]\text{ for }k\in{\Bbb N_0},          \text{
in particular,}\\
\Xi _+^\mu U_0&=\Cal F^{-1}[(\sigma +i\xi _n)^\mu \hat \varphi _0(\xi ')(\sigma +i\xi
_n)^{-\mu }]=\varphi _0(x')\otimes \delta (x_n),
\\
\Xi _+^\mu U_1&=\Cal F^{-1}[(\sigma +i\xi _n)^\mu \hat \varphi _1(\xi ')(\sigma +i\xi
_n)^{-\mu -1}]=\Cal F^{-1}_{\xi '\to x'}[\hat \varphi _1(\xi ')He^{-\sigma x_n}].
\endaligned\tag 6.10
$$

One more conclusion will be drawn; we here refer to the space $H^{(\mu
-1)(s)}(\crnp)$ recalled in Section 7, defined as  when $s>\mu -\frac32$.

\proclaim{Lemma 6.1} For $\mu >0$, $e^+x_n^{\mu -1}\Cal S(\crnp)= \Xi
_+^{-\mu +1}e^+\Cal S(\crnp)$.
\endproclaim

\demo{Proof} We have seen above that when $u\in e^+x_n^{\mu -1}\Cal
S(\crnp)$, then for any $k\in\N_0$, $u=U_0+\dots+U_k+u^{(k)}$, where
the $U_j$ are as defined in (6.7), and $u^{(k)}\in   e^+x_n^{\mu -1+k}\Cal
S(\crnp)$. Here
$$ 
U_j=\Cal F^{-1}[\hat\varphi _j(\sigma +i\xi _n)^{-\mu -j}]=\Xi
_+^{-\mu +1}\Cal F^{-1}(\hat\varphi _j(\sigma +i\xi _n)^{-j-1})\in \Xi
_+^{-\mu +1}e^+\Cal S(\crnp)
$$
 for all $j$. The remainder $
u^{(k)}$ is in $ e^+x_n^{\mu -1+k}\Cal S(\crnp)$, 
which for $k\to\infty $ converges to
$\dot {\Cal S}(\crnp)=\{v\in \Cal S(\rn)\mid \supp v\subset\crnp\}$,
where all boundary values at $x_n=0$ vanish. Here  $\Xi _+^{-\mu
+1}\dot {\Cal S}(\crnp)=\dot {\Cal S}(\crnp)$, since $\Xi _+^{-\mu
+1}$ is a homeomorphism on $\Cal S(\rn)$, and preserves support in
$\crnp$, and so does its inverse $\Xi _+^{\mu -1}$. Thus letting $k\to\infty $ in the
decomposition, we find the inclusion  $e^+x_n^{\mu -1}\Cal
S(\crnp)\subset
\Xi _+^{-\mu +1}e^+\Cal S(\crnp)$.

The opposite inclusion follows in a similar way: For $u$ given in $\Xi _
+^{-\mu +1}e^+\Cal S(\crnp)$, let $v\in e^+\Cal S(\crnp)$ be such that $u=\Xi _+^{-\mu
+1}v$. By a Taylor expansion in $x_n$ of $\Cal F_{x'\to \xi '}(e^{\sigma
x_n}v(x',x_n))$, we find an expansion of  $\acute v=\Cal F_{x'\to \xi
'}v$ in terms $I^j(x_n)\hat \varphi _j(\xi ')e^{-\sigma x_n}$, so that an application
of $\Xi_ +^{-\mu +1}$ to the corresponding term in $v$ gives $U_j=\Cal
F^{-1}[\hat \varphi _j(\sigma +i\xi _n)^{-\mu -j}]$. As we know from the
preceding analysis, these $U_j$ are also in
$e^+x_n^{\mu -1}\Cal S(\crnp)$. The remainder in $v$ after $k$ terms
is in $e^+x_n^k\Cal S(\crnp)$, which for $k\to\infty $ converges to
$\dot \Cal S(\crnp)$. $\Xi _+^{-\mu +1}$ can be applied along the way,
and in the limit we use that  $\Xi _+^{-\mu
+1}\dot {\Cal S}(\crnp)=\dot {\Cal S}(\crnp)$, to reach the
conclusion.
\qed
\enddemo

All the above is a straightforward use of the Fourier transform and
distribution theory. Further below,
we shall moreover
need to apply some basic constant-coefficient rules from the Boutet de
Monvel calculus. They are covered e.g.\ by the Appendix of \cite{G18}, which we shall
not repeat here.

Including that terminology, we observe that there are several
important descriptions of $U_0$ (recall that $\varphi _0=u_0$):
$$\aligned
U_0&=\Cal
F^{-1}_{\xi '\to x'}[\hat u_0(\xi ')I^{\mu -1}(x_n)e^{-\sigma
x_n}]
=\Cal F^{-1}_{\xi \to x}[\hat u_0(\xi ')(\sigma +i\xi
_n)^{-\mu }]\\
&=\Xi _+^{1-\mu }\Cal F^{-1}_{\xi \to x}[\hat u_0(\xi ')(\sigma +i\xi
_n)^{-1}]=\Xi _+^{1-\mu }e^+K_0u_0,\\
U_0&=I^{\mu -1}(x_n)\Cal
F^{-1}_{\xi '\to x'}[\hat u_0(\xi ')e^{-\sigma
x_n}]
=I^{\mu -1}(x_n)e^+K_0u_0= \tfrac1{\Gamma (\mu )}x_n^{\mu -1}e^+K_0u_0,
\endaligned\tag6.11
$$
where $K_0$ is the well-known Poisson operator
$
K_0\varphi =\Cal F^{-1}_{\xi '\to x'}[\hat \varphi (\xi ')r^+e^{-\ang{\xi '}
x_n}]
$,
defining a right-inverse of $\gamma _0$ satisfying $(1-\Delta )K_0=0$ on $\rnp$; its symbol is $(\ang{\xi '} +i\xi
_n)^{-1}$. There is also defined a related right-inverse
of $\gamma _0^{\mu -1}$, namely the operator
$$
K_0^{\mu -1}\equiv \Xi _+^{1-\mu }e^+K_0;\tag6.12
$$
it is a Poisson-like operator, solving the nonhomogeneous Dirichlet
problem for $(1-\Delta )^\mu $ with zero interior data. Moreover, in
view of (6.11),
$$
K_0^{\mu -1}=I^{\mu -1}(x_n)e^+K_0=\tfrac1{\Gamma (\mu )}x_n^{\mu -1}e^+K_0,
\tag6.13
$$
and $U_0=K_0^{\mu -1}u_0$.
There is a detailed study of the role of  $K_0^{\mu -1}$ and higher-order
Poisson-like operators in \cite{G19}.

\medskip
Now consider our operator $L$, written as $\Xi _-^{\mu ^*}Q\Xi _+^\mu $ as
in Section 4 (disregarding a smoothing term that does not contribute
to the formula). We shall show a Green's formula that allows
writing $\ang {r^+Lu,v}-\ang {u,r^+ L^*v}$ as a boundary integral when both $u$
and $v$
have nontrivial Dirichlet and Neumann traces. The boundary terms are
local, except for a term involving the Dirichlet traces of $u$ and $v$
with a pseudodifferential coefficient $(\mu -\mu ^*)\ang{D'}+B$.

\proclaim{Theorem 6.2} Let $L$,  $a$, $\delta
 $, $\mu $, $\mu ^*$ and $s_0$ be as described
 in Theorems {\rm 3.1} and {\rm 4.3}, and assume that
 $\mu ,\mu ^*>0$. For $u\in  \E_{\mu -1}(\crnp)\cap \E'(\rn)$,
$v\in  \E_{\mu ^*-1}(\crnp)\cap \E'(\rn)$,  there holds:
$$
\multline
\ang{r^+Lu,v}_{\ol H^{-\mu ^*+\frac12+\varepsilon }, \dot
H^{\mu ^*-\frac12-\varepsilon }}-\ang{u,r^+{L^*v}}_{\dot
H^{\mu -\frac12-\varepsilon },\ol H^{-\mu +\frac12+\varepsilon }}\\
=\ang{s_0\gamma _1^{\mu -1}u,\gamma _0^{\mu ^*-1}v
}
-\ang{s_0\gamma _0^{\mu -1}u ,\gamma
_1^{\mu ^*-1}v  }+\ang {[s_0(\mu -\mu ^*)\ang{D'}+B]\gamma _0^{\mu -1}u,\gamma _0^{\mu ^*-1}v},
\endmultline\tag6.14
$$
with $L_2({\Bbb R}^{n-1})$-scalar products in the right-hand side; it may
also be written
$$
\multline
\int_{\rnp} Lu\,\bar v\,dx-\int_{\rnp}u\,\overline{ L^*v}\,dx\\
=\int_{{\Bbb R}^{n-1}}
\bigl(s_0\Gamma (\mu +1)\Gamma (\mu ^*)\gamma _1(\tfrac {u}{x_n^{\mu -1 }})\,\gamma
_0(\tfrac{\bar v}{x_n^{\mu ^*-1}})-s_0\Gamma (\mu )\Gamma (\mu ^*+1)\gamma _0(\tfrac {u}{x_n^{\mu -1 }})\,\gamma
_1(\tfrac{\bar v}{x_n^{\mu ^*-1}})\\
+[s_0(\mu -\mu ^*)\ang{D'}+B]\Gamma (\mu )\Gamma (\mu ^*)\gamma _0(\tfrac u{x_n^{\mu -1 }})\,\gamma
_0(\tfrac{\bar v}{x_n^{\mu ^*-1}})\bigr)
\, dx'.
\endmultline
\tag6.15
$$
Here, with 
$$
q(\xi )=(\ang{\xi '}-i\xi _n)^{-\mu ^*}\Cal L(\xi )(\ang{\xi '}+i\xi
_n)^{-\mu },\tag 6.16
$$
the first-order pseudodifferential operator $B$ on ${\Bbb R}^{n-1}$
equals $\OP(b(\xi '))$, where $b(\xi ')$ is the jump at $x_n=0$ of $\Cal
F^{-1}_{\xi _n\to x_n}(q-s_0)$.  
\endproclaim

\example{Remark 6.3}
It is also shown in the proof, using terminology from the Boutet de Monvel
calculus, that when $q$ is decomposed as $q(\xi )=s_0+f_+(\xi
)+f_-(\xi )$ as in (4.14),
$f_+=h^+q$ and $f_-=h^-_{-1}q$, then
$$
b(\xi ')=\frak b(\xi ')-\overline{\frak b'}(\xi '), \text{ where } \frak b =\tfrac1{2\pi }\int^+ f_+\,
d\xi _n ,\; \frak b' =\tfrac1{2\pi }\int^+ \overline{f_-}\,
d\xi _n.\tag6.17
$$
\endexample

\demo{Proof of Theorem {\rm 6.2}}
Take $u\in  \E_{\mu -1}(\crnp)\cap \E'(\rn)$,
 $v\in  \E_{\mu ^*-1}(\crnp)\cap \E'(\rn)$, and write them
as 
$$
u=u'+U_0\text{ with }U_0=K^{\mu -1}_0u_0,\quad v=v'+V_0\text{ with }V_0=K_0^{\mu ^*-1}v_0,\tag6.18
$$
according to the above notation with $\mu $ for $u$, $\mu ^*$ for
$v$. Here $U_0=K^{\mu -1}_0u_0\in \Xi _+^{-\mu +1}e^+\Cal
S(\crnp)$\linebreak $=e^+x_n^{\mu -1}\Cal S(\crnp)$ (cf.\ Lemma 6.1),
whereas $u'$ lies in the smaller space $e^+x_n^\mu \Cal
S(\crnp)=  \Xi _+^{-\mu }e^+\Cal S(\crnp)$, as accounted for
before the theorem.

Since
$v\in \Xi _+^{1-\mu ^*}e^+\Cal S(\crnp)\subset \Xi _+^{1-\mu ^*}\dot H^{\frac12-\varepsilon }(\crnp)=\dot H^{\mu ^*-\frac12-\varepsilon }(\crnp)
$ and $r^+Lu\in C^\infty (\crnp) \cap \bigcap_s \ol H^s(\rnp)$, the
scalar product of $r^+Lu$ and $v$ can be interpreted as the duality
$$
I=\ang{r^+Lu,v}_{ \ol H^{\frac12- \mu ^*+\varepsilon }(\rnp) ,\dot H^{ \mu ^*-\frac12-\varepsilon
}(\crnp)},
$$
as in the proof of Theorem 4.1 above. It splits into four pieces:
$$
\aligned
I&=I_1+I_2+I_3+I_4, \\
I_1&=\ang{r^+Lu',v'},\quad I_2=\ang{r^+LK_0^{\mu -1}u_0
,v'},\\
I_3&=\ang{r^+Lu',K_0^{\mu ^*-1} v_0  },\quad I_4=\ang{r^+LK_0^{\mu -1} u_0 ,K_0^{\mu ^*-1} v_0  }.
\endaligned\tag6.19
$$

$I_1$ will be kept unchanged, to match a similar term with $L^*$ later.

For $I_2$ we observe, using (6.10), (6.11) and the representation
$Q=\OP(q)$, $q=s_0+f_++f_-$ from (4.14):
$$
\aligned
I_2&=\ang{r^+LK_0^{\mu -1} u_0 ,v'}_{\ol H^{-\mu ^*+\frac12+\varepsilon }, \dot
H^{\mu ^*-\frac12-\varepsilon }}\\
&=\ang{\Xi_{-,+}^{\mu ^*}r^+Q\Xi_+^{\mu }e^+K_{0}^{\mu -1}u_0
,v'}_{\ol H^{-\mu ^*+\frac12+\varepsilon }, \dot H^{\mu ^*-\frac12-\varepsilon
}}\\
&=\ang{r^+Q(u_0(x')\otimes\delta (x_n)) ,\Xi 
_+^{\mu ^*}v'}_{\ol H^{\frac12+\varepsilon }, \dot H^{-\frac12-\varepsilon
}}\\
&=\ang{r^+\operatorname{OPK}(f_+) u_0 ,\Xi 
_+^{\mu ^*}v'}_{\ol H^{\frac12+\varepsilon }, \dot H^{-\frac12-\varepsilon }}.
\endaligned\tag6.20
$$
The last equality came from the rule from the Boutet de Monvel calculus that the
mapping $u_0\mapsto r^+Q(u_0(x')\otimes\delta (x_n))  $ is the Poisson
operator with symbol $h^+q(\xi )=f_+(\xi )$; the contributions from $s_0$ and
$f_-(\xi )$ vanish since they map into distributions supported in $\crnm$.

Next, consider $I_3$. Here, by (6.12) and (4.2), 
$$
\aligned
I_3&=\ang{r^+Lu',K_0^{\mu ^*-1} v_0  }=\ang{\Xi_{-,+}^{\mu ^*}r^+Q\Xi_+^\mu u',\Xi_+^{1-\mu ^*}e^+K_{0} v_0  }_{\ol H^{-\mu ^*+\frac12+\varepsilon }, \dot
H^{\mu ^*-\frac12-\varepsilon }}\\
&=\ang{\Xi_{-,+}^{1-\mu ^*}\Xi _{-,+}^{\mu ^*}r^+Q\Xi_+^\mu u',K_{0} v_0  }_{\ol H^{-\frac12+\varepsilon }, \dot
H^{\frac12-\varepsilon }}\\
&=\ang{\Xi_{-,+}^1r^+Q\Xi_+^\mu u',K_{0} v_0  }_{\ol H^{-\frac12+\varepsilon }, \dot
H^{\frac12-\varepsilon }}\\
&=\ang{K_0^*\Xi_{-,+}^1Q\Xi_+^\mu u', v_0  }_{ H^{\varepsilon }({\Bbb R}^{n-1}), 
H^{-\varepsilon }({\Bbb R}^{n-1})}
=(K_0^*\Xi_{-,+}^1Q\Xi_+^\mu u', v_0  )_{L_2({\Bbb R}^{n-1})},
\endaligned\tag6.21
$$
since $ v_0  $ is in $L_2({\Bbb R}^{n-1})$.
It is
used that $\ol H^t$ identifies with $\dot H^t$ for $|t|<\frac12$
(there the
indication $e^+$ is understood). 

Denote $\Xi_+^\mu u'=w\in e^+\ol H^{\frac12+\varepsilon
}(\rnp)$. Since $K_0^*$ is the trace operator with symbol \linebreak
 $(\ang{\xi'}-i\xi _n)^{-1}$, the rules of calculus give that
$$
\multline
K_0^*r^+\Xi_-^1Qw=\operatorname{OPT}(h^-((\ang{\xi'}-i\xi
_n)^{-1}(\ang{\xi'}-i\xi _n)q(\xi ))w\\
=\operatorname{OPT}(h^-q(\xi ))w=\operatorname{OPT}(s_0+f_-(\xi ))w,
\endmultline
$$
and hence, in view of (6.9), (6.10), 
$$
\aligned
I_3&=(\operatorname{OPT}(s_0+f_-)\Xi_+^\mu u', v_0 
)_{L_2({\Bbb R}^{n-1})}=((s_0\gamma _0+\operatorname{OPT}(f_-))\Xi_+^\mu u', v_0 )_{L_2({\Bbb R}^{n-1})}\\
&=
(s_0\gamma _0^\mu u', v_0 )_{L_2({\Bbb
R}^{n-1})}+(\operatorname{OPT}(f_-)\Xi_+^\mu u', v_0  )_{L_2({\Bbb
R}^{n-1})}\\
&=
(s_0(\gamma _1^{\mu -1}u+\mu \ang{D'}\gamma _0^{\mu -1}u, v_0 )_{L_2({\Bbb R}^{n-1})}+(\operatorname{OPT}(f_-)\Xi_+^\mu u', v_0  )_{L_2({\Bbb R}^{n-1})}.
\endaligned\tag6.22
$$

Finally, consider $I_4$. Beginning as in the treatment of $I_3$, we find:
$$
\aligned
 I_4&=\ang{r^+LK_0^{\mu -1} u_0 ,K_0^{\mu ^*-1} v_0  }_{\ol H^{-\mu ^*+\frac12+\varepsilon }, \dot
H^{\mu ^*-\frac12-\varepsilon }}\\
&=\ang{\Xi 
 _{-,+}^{\mu ^*}r^+Q\Xi_+^\mu \Xi_+^{1-\mu }e^+K_0 u_0 ,\Xi 
 _+^{1-\mu ^*}e^+K_0 v_0  }_{\ol H^{-\mu ^*+\frac12+\varepsilon }, \dot
H^{\mu ^*-\frac12-\varepsilon }}\\
&=\ang{\Xi_{-,+}^1r^+Q\Xi_+^1e^+K_0 u_0 ,K_0 v_0  }_{\ol H^{-\frac12+\varepsilon }, \dot
H^{\frac12-\varepsilon }}\\
&=\ang{K_0^*\Xi_{-,+}^1r^+Q\Xi_+^1e^+K_0 u_0 , v_0  }_{ H^{\varepsilon }({\Bbb R}^{n-1}), 
H^{-\varepsilon }({\Bbb R}^{n-1})}\\
&=(\Cal B u_0 , v_0  )_{L_2({\Bbb R}^{n-1})},
\endaligned\tag6.23
$$
where ${\Cal B}=K_0^*\Xi_{-,+}^1r^+Q\Xi_+^1e^+K_0$ is a certain $\psi $do
on ${\Bbb R}^{n-1}$ of order 1. We can reduce this expression by rules
of calculus involving the so-called  plus-integral, cf.\ \cite{G18,(A.14)ff., (A.15)}. 
Since $K_0$ has symbol $(\ang{\xi '}+i\xi _n)^{-1}$, the symbol of the Poisson operator
$r^+Q\Xi _+^{1}e^+K_0$ is $h^+q=f_+$, which by composition 
with $\Xi_{-,+}^1$ to the left gives a Poisson operator with
symbol $h^+((\ang{\xi '}-i\xi _n)f_+)$. The symbol $\frak b(\xi ')$ of ${\Cal B}$  is calculated as in
\cite{G18,(4.15)}:
$$
\aligned
\frak b(\xi ')&=\tfrac1{2\pi }\int^+(\ang{\xi '}-i\xi _n)^{-1}h^+((\ang{\xi '}-i\xi _n)f_+)\,
d\xi _n\\
&=\tfrac1{2\pi }\int^+ f_+\,
d\xi _n = \lim_{z_n\to
0+}\check {q}(\xi ',z_n),
\endaligned\tag6.24
$$
where  $\check {q}(\xi ',z_n)$ stands for
 $\Cal F^{-1}_{\xi _n\to z_n}q(\xi )$.

The same arguments apply to $I'=\ang{u,r^+L^*v}$, after a conjugation and
an exchange of $\mu ,\mu ^*$ by $\mu ^*,\mu $. This gives:
$$
\aligned
I'&=I'_1+I'_2+I'_3+I'_4,\text{ where }\\
I'_1&=\ang{u',r^+P^*v'}_{\dot
H^{\mu -\frac12-\varepsilon },\ol H^{-\mu +\frac12+\varepsilon }},\\
I'_2&=\ang{\Xi 
_+^\mu u',\operatorname{OPK}(\overline{f_-}) v_0 , }_{\dot
H^{-\frac12-\varepsilon },\ol H^{\frac12+\varepsilon }},\\
I_3'&=\ang{s_0 u_0 ,\gamma _1^{\mu ^*-1}v+\mu ^*\ang{D'}\gamma _0^{\mu ^*-1}v}
+\ang{ u_0  ,\operatorname{OPT}(\overline{f_+})\Xi_+^{\mu ^*}v' },\\
I_4'&=( u_0 ,{\Cal B}' v_0  )_{L_2({\Bbb R}^{n-1})},
\endaligned\tag6.25
$$
where $\Cal B'=\OP(\frak b'(\xi ))$, with
$$
\frak b'(\xi ')
=\tfrac1{2\pi }\int^+h^+\bar q(\xi )\,
d\xi _n =\tfrac1{2\pi }\int^+\overline{f_-}\,
d\xi _n ,
\quad \overline{\frak b'}= \lim_{z_n\to
0-}\check q(\xi ',z_n).
\tag6.26
$$

 Now to calculate $I-I'$, we first find
$$
I_1-I'_1=0
$$
by Theorem 5.3, which clearly covers $u'\in e^+x_n^\mu \Cal S(\crnp)$, $v'\in e^+x_n^{\mu ^*}\Cal S(\crnp)$.  Next,
$$
\aligned
I_2-I'_3&=\ang{\operatorname{OPK}(f_+)u_0, \Xi
_+^{\mu ^*}v'}
-\ang{s_0 u_0 ,\gamma _1^{\mu ^*-1}v+\mu ^*\ang{D'}\gamma _0^{\mu
^*-1}v}
-\ang{ u_0 ,\operatorname{OPT}(\overline{f_+})\Xi_+^{\mu ^*}v' }
\\
&=\ang{\operatorname{OPK}(f_+)u_0, \Xi
_+^{\mu ^*}v'}-\ang{\operatorname{OPK}({f_+}) u_0 ,\Xi_+^{\mu ^*}v' }
- \ang{s_0 u_0 ,\gamma _1^{\mu ^*-1}v+\mu ^*\ang{D'}\gamma _0^{\mu
^*-1}v}\\
&=-\ang{s_0 u_0 ,\gamma _1^{\mu ^*-1}v+\mu ^*\ang{D'}\gamma _0^{\mu ^*-1}v},
\endaligned
$$
using that the adjoint of the trace operator
$\operatorname{OPT}(\overline{f_+})$ is the Poisson operator $\operatorname{OPK}(f_+)$. There is a similar calculation of $I_3-I'_2$, so we find
$$
I_2+I_3-I'_2-I'_3=\ang{s_0\gamma _1^{\mu -1}u,\gamma _0^{\mu ^*-1}v
}-\ang{s_0\gamma _0^{\mu -1}u ,\gamma
_1^{\mu ^*-1}v  }+\ang {s_0(\mu -\mu ^*)\ang{D'}\gamma _0^{\mu -1}u,\gamma _0^{\mu ^*-1}v}.
$$
(Since $u'$ and $v'$ are solutions of homogeneous Dirichlet problems,
whereas $K_0^{\mu -1}u_0$ and $K_0^{\mu ^*-1}v_0$ solve nonhomogeneous
Dirichlet problems, this could also have been derived using the halfways Green's formula.)
Finally,
$$
I_4-I'_4=(({\Cal B}-{{\Cal B}'}^*)u_0,v_0)_{L_2({\Bbb R}^{n-1})},=(({\Cal B}-{{\Cal B}'}^*)\gamma _0^{a-1}u,\gamma _0^{a-1}v)_{L_2({\Bbb R}^{n-1})},
$$
where  $B={\Cal B}-{{\Cal B}'}^*$ satisfies, with $b(\xi ')=\frak
b(\xi ')-\overline{\frak b'}(\xi ')$,
$$
B=\OP(b(\xi ')),\quad b(\xi ')=\lim_{z_n\to
0+}\check q(\xi ',z_n)-\lim_{z_n\to
0-}\check q(\xi ',z_n),\tag6.27
$$
the jump at $z_n=0$ in the bounded part of $\check q(\xi ',z_n)$.
Adding the terms, we find the assertion in the theorem. \qed

 \enddemo

\subhead 7. Results in other function spaces \endsubhead

For the interested reader, we shall go rapidly through some results in
function spaces of finite regularity, that hold along with the above
statements shown for $C^\infty $-functions.

 First we recall the definitions of the spaces; more details are given in
 \cite{G15,G14}.

The $\mu $-transmission spaces $H^{\mu (s)}(\crnp)$ are defined as $\Xi _+^{-\mu }e^+\ol
H^{s-\mu }(\rnp)$ for $s>\mu -\frac12$, and satisfy when $\mu
>-1$:
$$
H^{\mu (s)}(\crnp)\cases =\dot H^s(\crnp) \text{ if }s\in
\,]\mu -\frac12,\mu +\frac12[\,,\\
\subset e^+x_n^\mu  \ol H^{s- \mu }(\rnp)+\dot
H^{s(-\varepsilon )}(\crnp)\text{ if }s>\mu +\frac12, 
\endcases
\tag 7.1
$$
 where $-\varepsilon $ is active if $s-\mu -\frac12$ is integer. Here
 the trace operator $\gamma _k^\mu $ maps $H^{\mu (s)}(\crnp)$
 continuously into  $ H^{s-\mu -\frac12}({\Bbb R}^{n-1})$ when $s>\mu +k
 +\frac12$; note that the range is in $L_2(\R^{n-1})$.

\proclaim{Lemma 7.1} $\Cal E_{\mu }(\crnp)\cap \E'(\rn)$ is dense in
$H^{\mu (s)}(\crnp)$, for all $s>\mu -\frac12$.
\endproclaim

\demo{Proof} This is a corollary to Lemma 6.1, where the identity $e^+x_n^{\mu }\Cal S(\crnp)= \Xi
_+^{-\mu }e^+\Cal S(\crnp)$ for $\mu >-1$ was proved. Note that since $\Cal
S(\crnp)$ is densely embedded in $\ol H^{s-\mu }(\rnp)$, $\Xi _+^{-\mu }e^+\Cal
S(\crnp)=e^+x_n^{\mu }\Cal S(\crnp)$ is densely embedded in $\Xi
_+^{-\mu }e^+\ol H^{s-\mu}(\rnp)=H^{\mu (s)}(\crnp)$. When $\varphi
\in C_0^\infty (\rn)$ with
$\varphi (x)=1$ for $|x|\le 1$, then if $v\in e^+x_n^{\mu }\Cal S(\crnp)$, $\varphi (\delta x)v\to v$ in
$e^+x_n^{\mu }\Cal S(\crnp)$ for $\delta \to 0$, and a fortiori in
$ H^{\mu (s)}(\crnp)$. For a $u\in H^{\mu (s)}(\crnp)$, one
can for any $k\in{\Bbb N}$ find
$u_k\in e^+x_n^{\mu }\Cal S(\crnp)$ such that
$\|u-u_k\|_{H^{\mu (s)}}\le 1/k$; then $\|u-\varphi (\delta x)u_k\|_{H^{\mu (s)}}\le 2/k$ for $\delta $ sufficiently small, and $\varphi (\delta
x)u_k\in \E_{\mu -1}(\crnp)\cap \E'(\rn)$.
\qed
\enddemo

\proclaim{Lemma 7.2} The mapping $\partial_n$ sends $H^{\mu(s
)}(\crnp)$ continuously into $ H^{(\mu
-1)(s-1)}(\crnp)$ for all $\mu $, all $s>\mu +\frac12$.
\endproclaim

\demo{Proof}
This follows from
$$
\aligned
\partial_nH^{\mu(s)}(\crnp)&=\partial_n\Xi _+^{-\mu }e^+\ol H^{s-\mu
}(\rnp)=
(\partial_n+\ang{D'}-\ang{D'})\Xi _+^{-\mu }e^+\ol H^{s-\mu }(\rnp)\\
&\subset \Xi _+^{-\mu +1}e^+\ol H^{s-\mu }(\rnp)+\Xi _+^{-\mu }e^+\ol
H^{s-1-\mu }(\rnp)\\
&=H^{(\mu -1)(s-1)}(\crnp)+H^{(\mu
)(s-1)}(\crnp)=H^{(\mu -1)(s-1)}(\crnp),
\endaligned
$$
using that $\partial_n+\ang{D'}=\OP(i\xi _n+\ang{\xi '})=\Xi _+^1$. \qed
\enddemo

\proclaim{Theorem 7.3} For $s>\mu +\frac12$, $s'>\mu ^*+\frac12$,
Theorem {\rm 4.3} extends to $u\in H^{\mu (s)}(\crnp)$, 
$u'\in H^{\mu ^*(s')}(\crnp)$, 
in the form 
$$
\aligned
\ang{r^+Lu,\partial_nu'}&_{\ol H^{-\mu
^*+\frac12+\varepsilon}(\rnp),\dot H^{\mu ^*-\frac12-\varepsilon}(\crnp)}+\ang{\partial_nu,L^*u'}_{\dot H^{\mu -\frac12-\varepsilon}(\crnp),\ol H^{-\mu
+\frac12+\varepsilon}(\rnp)}\\
&=\Gamma (\mu +1){\Gamma(\mu ^*+1)}\int_{{\Bbb
R}^{n-1}}s_0\gamma _0(u/x_n^{\mu })\,{\gamma _0(\bar u'/x_n^{\mu ^*})}\,
dx'.
\endaligned\tag7.2
$$
\endproclaim

\demo{Proof} When $u'\in H^{\mu ^*(s')}(\crnp)$ with  $s'>\mu
^*+\frac12$, then since (a fortiori)  $s'-\mu ^*>\frac12-\varepsilon $, we have using Lemma 7.2:
$$
\aligned
\partial_nu'\in  H^{(\mu ^*-1)(s'-1)}(\crnp)&=\Xi _+^{-\mu ^*+1}e^+\ol
H^{s'-\mu ^*}(\crnp)\subset \Xi _+^{-\mu ^*+1}e^+\ol
H^{\frac12-\varepsilon }(\rnp)\\
&=\Xi _+^{-\mu ^*+1}\dot
H^{\frac12-\varepsilon }(\crnp)=\dot H^{\mu ^*-\frac12-\varepsilon }(\crnp),
\endaligned
$$
 for small $\varepsilon >0$. At the same time, when $u\in H^{\mu
 (s)}(\crnp)$ with $s>\mu
+\frac12$, then (by \cite{G15, Th.\ 4.2})
$$
r^+Lu\in \ol H^{s-2a }(\rnp)\subset \ol H^{\mu +\frac12-2a+\varepsilon
}(\rnp)
=\ol H^{\,-\mu ^*+\frac12+\varepsilon }(\rnp)
$$
 for small $\varepsilon >0$. Thus the first duality in (7.2) has a
 sense for such $u,u'$.

Similarly, the second duality in (7.2) has a sense when $s>\mu
+\frac12$, $s'>\mu ^*+\frac12$.

In view of Lemma 7.1,
we can choose sequences of
approximating functions $(u_k)_{k\in\N}\subset \E_\mu(\crnp)\cap \E'(\rn) $ converging to $u$, $(u'_k)_{k\in\N}\subset
\E_{\mu^*}\crnp)\cap \E'(\rn) $ converging to $u'$, in the
abovementioned norms.
The  formula (7.2)
holds for the pair $\{u_k,u'_k\}$ by Theorem 4.3, each $k$, and converges to the formula for the given
$u,u'$ when $k\to\infty $. \qed
\enddemo

The halfways Green's formula similarly extends:

\proclaim{Theorem 7.4} Let $\mu ,\mu ^*>0$.  For  $s>\mu +\frac12$,
$s'>\mu ^*+\frac12$, Theorem {\rm 5.1} extends to $u\in H^{(\mu
-1)(s)}(\crnp)$, 
$v\in H^{\mu ^*(s')}(\crnp)$, 
in the form 
$$
\aligned
\ang{r^+Lu,v}&_{\ol H^{-\mu
^*-\frac12+\varepsilon}(\rnp),\dot H^{\mu ^*+\frac12-\varepsilon}(\crnp)}-\ang{u,L^*v}_{\dot H^{\mu -\frac12-\varepsilon}(\crnp),\ol H^{-\mu
+\frac12+\varepsilon}(\rnp)}\\
&=-\Gamma (\mu ){\Gamma(\mu ^*+1)}\int_{{\Bbb
R}^{n-1}}s_0\gamma _0(u/x_n^{\mu -1 })\,{\gamma _0(\bar v/x_n^{\mu ^*})}\,
dx'.
\endaligned
\tag7.3
$$
\endproclaim

\demo{Proof}  When $u\in H^{(\mu -1)(s)}(\crnp)$ with  $s>\mu
+\frac12$, then since (a fortiori)  $s-\mu +1 >\frac12-\varepsilon $, we have:
$$
\aligned
u\in  H^{(\mu -1)(s)}(\crnp)&=\Xi _+^{-\mu +1}e^+\ol
H^{s-\mu +1}(\crnp)\subset \Xi _+^{-\mu +1}e^+\ol
H^{\frac12-\varepsilon }(\rnp)\\
&=\Xi _+^{-\mu +1}\dot
H^{\frac12-\varepsilon }(\crnp)=\dot H^{\mu -\frac12-\varepsilon }(\crnp),
\endaligned
$$
 for small $\varepsilon >0$. At the same time, when $v\in H^{\mu ^*
  (s')}(\crnp)$ with $s'>\mu ^*
+\frac12$, then 
$$
r^+L^*u\in \ol H^{s'-2a }(\rnp)\subset \ol H^{\mu ^*+\frac12-2a+\varepsilon
}(\rnp)
=\ol H^{\,-\mu +\frac12+\varepsilon }(\rnp)
$$
 for small $\varepsilon >0$. Thus the second duality in (7.3) has a
 sense for such $u,v$.

We also have, since $v\in H^{\mu ^*
  (s')}(\crnp)$ with $s'>\mu ^*
+\frac12$, that
$$
\aligned
v\in  H^{\mu ^*(s')}(\crnp)&=\Xi _+^{-\mu ^*}e^+\ol
H^{s'-\mu ^*}(\crnp)\subset \Xi _+^{-\mu ^*}e^+\ol
H^{\frac12-\varepsilon }(\rnp)\\
&=\Xi _+^{-\mu ^*}\dot
H^{\frac12-\varepsilon }(\crnp)=\dot H^{\mu ^*+\frac12-\varepsilon }(\crnp),
\endaligned
$$
and since  $u\in H^{(\mu -1)
 (s)}(\crnp)$ with $s>\mu -1
+\frac12=\mu -\frac12$, then (by \cite{G15, Th.\ 4.2})
$$
r^+Lu\in \ol H^{s-2a }(\rnp)\subset \ol H^{\mu -\frac12-2a+\varepsilon
}(\rnp)
=\ol H^{\,-\mu ^*-\frac12+\varepsilon }(\rnp)
$$ 
 for small $\varepsilon >0$. Thus the first duality in (7.3) is
 well-defined.

The formula (7.3)  is then deduced from (5.1) by approximation, as in
the proof of Theorem 7.3. \qed
\enddemo

By a similar proof, the full Green's formula extends:

\proclaim{Theorem 7.5} Let $\mu ,\mu ^*>0$.  For  $s>\mu +\frac12$,
$s'>\mu ^*+\frac12$, Theorem {\rm 6.2} extends to $u\in H^{(\mu
-1)(s)}(\crnp)$, 
$v\in H^{(\mu ^*-1)(s')}(\crnp)$, 
in the form
$$
\multline
\ang{r^+Lu,v}_{\ol H^{-\mu ^*+\frac12+\varepsilon }, \dot
H^{\mu ^*-\frac12-\varepsilon }}-\ang{u,r^+{L^*v}}_{\dot
H^{\mu -\frac12-\varepsilon },\ol H^{-\mu +\frac12+\varepsilon }}\\
=\ang{s_0\gamma _1^{\mu -1}u,\gamma _0^{\mu ^*-1}v
}
-\ang{s_0\gamma _0^{\mu -1}u ,\gamma
_1^{\mu ^*-1}v  }+\ang {[s_0(\mu -\mu ^*)\ang{D'}+B]\gamma _0^{\mu -1}u,\gamma _0^{\mu ^*-1}v},
\endmultline\tag7.4
$$
with $L_2({\Bbb R}^{n-1})$-scalar products in the right-hand side.
\endproclaim

As a corollary, it is found that the formulas hold for $u,u',v$ in
certain  H\"older-type spaces with parameter $s$,
contained in the Sobolev-type spaces with parameter $s-\varepsilon $. The transmission
spaces  are defined in \cite{G14} by the usual formula $C_*^{\mu (s)}(\crnp)=\Xi
_+^{-\mu }e^+\ol C_*^{s-\mu }(\rnp)$ for $s>\mu -1$, where $C_*^s$ stands for the standard
H\"older space when $s\in\rp\setminus \N$ and the H\"older-Zygmund
space otherwise, and
there is a description like (7.1):
$$
C_*^{\mu (s)}(\crnp)\cases =\dot C_*^s(\crnp) \text{ if }s\in
\,]\mu -1,\mu [\,,\\
\subset e^+x_n^\mu  \ol C_*^{s- \mu }(\rnp)+\dot
C_*^{s(-\varepsilon )}(\crnp)\text{ if }s>\mu , 
\endcases
\tag 7.5
$$
 where $-\varepsilon $ is active if $s-\mu $ is integer. For $s>\mu
 +k$, $k\in\N_0$, $\gamma _k^{\mu }$ maps $C_*^{\mu (s)}(\crnp)$ into
 $C_{*}^{s-\mu -k}(\R^{n-1})$. 
  The operator
$r^+L$ maps $C_*^{\mu (s)}(\crnp)$ into $\ol C_*^{s-2a}(\rnp)$ when
$s>\mu -1$. When $s>\mu >0$ and $s,s-\mu \notin\N$, $u\in C_*^{\mu
(s)}(\crnp)\implies u/x_n^\mu \in \ol C^{s-\mu }(\rnp)$ in view of (7.5).

Since $C^s_*(\rn)\subset H^{s-\varepsilon
}(\rn)$ (any $\varepsilon >0$) with
consequential embeddings for all the derived spaces, we can formulate, as special cases of Theorems 7.3--7.5:

\proclaim{Corollary 7.6} Let  $s>\mu +\frac12$ and $s'>\mu
^*+\frac12$, $s$ and $s' > 2a$.

 The formula in Theorem {\rm 7.3} holds when
$u\in C_*^{\mu (s)}(\crnp)$, 
$u'\in C_*^{\mu ^*
(s')}(\crnp)$. 

When $\mu ,\mu ^*>0$, the formula in Theorem {\rm 7.4} holds when
$u\in C_*^{(\mu -1)(s)}(\crnp)$, 
$v\in C_*^{\mu ^*(s')}(\crnp)$, 
and the formula in Theorem {\rm 7.5} holds when
$u\in C_*^{(\mu -1)(s)}(\crnp)$, 
$v\in$ \linebreak$ C_*^{(\mu ^*-1)
(s')}(\crnp)$. 
\endproclaim

Since $s,s'>2a$ here, $L$ and $L^*$ map into the continuous
functions on $\crnp$, so the formulas can be written with integrals
over $\rnp$ in the place of dualities, as
in Theorem {\rm 4.3} (when $\mu ,\mu ^*>0$), Theorem {\rm 5.1} and
Theorem {\rm 6.2}.


\enddocument